\documentclass[11pt]{article}
\usepackage{amsfonts,latexsym,amsmath,amscd,geometry}
\usepackage{amssymb}
\usepackage{latexsym}
\usepackage{xcolor}
\usepackage[title]{appendix}

\newcommand \nc{\newcommand}
\newtheorem{theorem}{Theorem}[section]
\newtheorem{lemma}[theorem]{Lemma}
\newtheorem{proposition}[theorem]{Proposition}

\newtheorem{remark}[theorem]{Remark}

\nc{\ba}{\begin{array}}\nc{\ea}{\end{array}}
\nc{\be}{\begin{eqnarray}}\nc{\ee}{\end{eqnarray}}
\nc{\beq}{\begin{equation}}\nc{\eeq}{\end{equation}}
\nc{\bex}{\begin{eqnarray*}}\nc{\eex}{\end{eqnarray*}}
\nc{\btm}{\begin{theorem}} \nc{\etm}{\end{theorem}}
\nc{\blm}{\begin{lemma}} \nc{\elm}{\end{lemma}}
\nc{\R}{\mathbb{R}}

\def\pf{\noindent{\bf Proof.\quad}}\def\endpf{\hfill$\Box$}
\def\di{\mbox{div\,}}
\def\curl{\mbox{curl\,}}

\def\u{\mathbf{u}}
\def\d{\mathbf{d}}
\def\la{\lambda}
\def\v{\varepsilon}
\begin{document}

\title{Singularity formation for full Ericksen-Leslie system of nematic liquid crystal flows in dimension two}

\author{Geng Chen
\footnote{Department of Mathematics, University of Kansas, Lawrence, KS 66045, U.S.A. Email: gengchen@ku.edu}\quad
Tao Huang
\footnote {Department of Mathematics, Wayne State University, Detroit, MI, 48202, U.S.A. Email: taohuang@wayne.edu}\quad 
Xiang Xu
\footnote{Department of Mathematics, University of Kansas, Lawrence, KS 66045, U.S.A. Email: x354x351@ku.edu}
\date{}
}
\maketitle

\begin{abstract} 
In this paper, we prove the singularity formation for Poiseuille laminar flow of full Ericksen-Leslie system modeling nematic liquid crystal flows in dimension two. The singularity is due to the geometric effect at the origin.
\end{abstract}

\section{Introduction}
Liquid crystal materials  have a degree of crystal structures but also  exhibit many hydrodynamic features  so they are capable to flow.  
Nematic liquid crystals are composed of rod-like molecules characterized by average alignment of  the long axes of neighboring molecules, which have the simplest structure  among liquid crystals and have been widely studied analytically and experimentally that lead to fruitful applications \cite{DeGP, Cha, Eri76, Les}. 

The dynamic theory of nematic liquid crystals  was first proposed by Ericksen \cite{ericksen62} and Leslie \cite{leslie68}  in the 1960s, (see, e.g. \cite{Les, lin89}).  
More precisely, if the orientation order parameters of molecules are treated as unit vector ${\d(x)}$, the classic static theory relies on the following Oseen-Frank free energy density
$$
2W(d,\nabla \d)
=k_1(\di \d)^2+k_2(\d\cdot\curl \d)^2+k_3|\d\times\curl \d|^2+(k_2+k_4)[\mbox{tr}(\nabla \d)^2-(\di \d)^2 ]
$$
where $k_i$, $i=1,\cdots, 3$ are the positive constants representing splay, twist, and bend affects respectively, with
$
k_2\geq |k_4|, \quad 2k_1\geq k_2+k_4.
$
The full Ericksen-Leslie system is given as follows

\begin{equation}\label{wels}
\begin{cases}
\dot \u+\nabla P=\nabla\cdot\sigma-\nabla\cdot\left(\frac{\partial W}{\partial\nabla \d}\otimes\nabla \d\right)
\\
\nabla\cdot \u=0\\
\nu \ddot{\d}=\gamma \d- g-\frac{\partial W}{\partial  \d}+\nabla\cdot\left(\frac{\partial W}{\partial\nabla \d}\right), \\
|\d|=1
\end{cases}
\end{equation}
where $\u$ is the velocity field of underlying incompressible fluid,
$\dot{f}=f_t+\u\cdot \nabla f$ is the material derivative,
$\gamma$ is Lagrangian multiplier of the constraint $|\d|=1$, and the constant $\nu>0$.
Let
$$
D= \frac12(\nabla \u+(\nabla \u)^T),\quad \omega= \frac12(\nabla \u-(\nabla \u)^T)=\frac12(\partial_j\u^i-\partial_i\u^j),\quad N=\dot \d-\omega \d,
$$
represent the rate of strain tensor, skew-symmetric part of the strain rate, and the
rigid rotation part of director changing rate by fluid vorticity, respectively. The kinematic transport $g$ is given by
$$
g=\lambda_1 N +\lambda_2Dd
$$
which represents the effect of the macroscopic flow field on the microscopic structure. The material coefficients $\lambda_1$ and $\lambda_2$ reflect the molecular shape and the slippery part between fluid and particles. The first term of $g$ represents the rigid rotation of molecules, while the second term stands for the stretching of molecules by the flow.
The viscous (Leslie) stress tensor $\sigma$ has the following form
$$
\sigma= \mu_1 (\d^TD\d)\d\otimes \d + \mu_2N\otimes \d   + \mu_3\d\otimes N+ \mu_4D + \mu_5(D\d)\otimes \d+ \mu_6\d\otimes (D\d) .
$$
These coefficients $\mu_i (1 \leq i \leq 6)$, depending on material and temperature, are called Leslie coefficients, and are related to certain local correlations in the fluid. The following relations are assumed in the literature.
$$
\lambda_1 =\mu_3-\mu_2,\quad \lambda_2 =\mu_6 -\mu_5,\quad \mu_2+ \mu_3 =\mu_6-\mu_5.
$$
The first two relations are compatibility conditions, while the third relation is called Parodi's relation, derived from Onsager reciprocal relations expressing the equality of certain relations between flows and forces in thermodynamic systems out of equilibrium. They also satisfy the following empirical relations (p.13, \cite{Les}) 
\begin{align}\label{alphas}
&\mu_4>0,\quad 2\mu_1+3\mu_4+2\mu_5+2\mu_6>0,\quad \gamma_1=\mu_3-\mu_2>0,\\
&  2\mu_4+\mu_5+\mu_6>0,\quad 4\gamma_1(2\mu_4+\mu_5+\mu_6)>(\mu_2+\mu_3+\gamma_2)^2\notag.
\end{align}
Note that the 4th relation  is implied by the 3rd together with the last relation.


The full Ericksen-Leslie system \eqref{wels} is a coupled system of forced Navier-Stokes equations and the wave map equations. Basic questions about existence, uniqueness and regularity of solutions are not completely understood. Several wellposedness results have been established in \cite{cai-wang20, HJLZ21, jiangluo17}, where local-in-time existence and uniqueness for initial data with finite energy, and global existence and uniqueness of classical solutions with small initial data have been established for various choices of coefficients. 

When $\nu=0$, the  Ericksen-Leslie system \eqref{wels} becomes a parabolic system (also called Ericksen-Leslie system in literature).  For the parabolic Ericksen-Leslie system in dimension two, the existence and uniqueness of global solution have been studied in \cite{WZZ13, huanglinwang14, hongxin12, LTX16, wangwang14}.   In dimension three, under some simplified assumptions, the authors of \cite{WZZ13} established global existence of solutions for small initial data and provided a characterization of the maximal existence time for general initial data.  We also refer the reader to the pioneer paper by Fanghua Lin  \cite{lin89}.

However, the full Ericksen-Leslie system \eqref{wels}, which is a coupled system of  forced Navier-Stokes equations and the wave map equations, is still poorly understood. In general, finite time singularity might form for solutions of \eqref{wels} even in one space dimension with smooth initial data. This new phenomenon on formation of cusp singularity due to the quasi-linearity in the wave equation, was first found by the authors in \cite{CHL20} for the one dimensional Poiseuille laminar flow. On the other hand, the global existence for H\"older continuous solution beyond cusp singularity was established in \cite{CHL20}. Also see \cite{CLS}. To obtain these results, we use earlier frameworks for the variational wave equation: formation of cusp singularity \cite{GHZ,CHL,CZ12}; existence \cite{BZ,HR,CZ12,BH,ZZ03}. Also see uniqueness \cite{BCZ}, and Lipschitz continuous dependence results in \cite{BC2015,BHY,BC}. 

Currently, the study of the solution in multiple space dimension is widely open. In this paper, we investigate the axisymmetric solutions to the full Ericksen-Leslie system \eqref{wels} with the special case of Oseen-Frank energy density
$$
2W(\d,\nabla \d)=|\nabla \d|^2,
$$
and a simple case of the coefficients  
\beq\label{specialmus}
\nu=1,\quad \mu_1=\mu_5=\mu_6=0,\quad \mu_2=-1,\quad \mu_3=1,\quad \mu_4=1, \quad \lambda_1=\mu_3-\mu_2=2.
\eeq
By the Onsager-Parodi relation,
\[
\lambda_2=\mu_6-\mu_5=\mu_2+\mu_3=0.
\]
The system \eqref{wels} becomes
\begin{equation}\label{fels}
\begin{cases}
 \u_t+\u\cdot \nabla \u+\nabla P=\nabla\cdot\big( \d\otimes N-N\otimes \d   + D\big)-\nabla\cdot\left(\nabla \d\odot\nabla \d\right)
\\
\nabla\cdot \u=0\\
\ddot{\d}+2\dot \d-2\omega \d=\Delta \d+\big(|\nabla \d|^2-|\dot{\d}|^2\big)\d. \\
|\d|=1.
\end{cases}
\end{equation}

The goal of this paper is to study the large data solution for \eqref{fels}. Since the  ${\bf d}$ equation has a ``semi-linear'' structure instead of a quasilinear one, we guess the finite time cusp singularity found in \cite{CHL20} will not form. However, in this paper, we prove that the systems will have another type of finite time blowup, similar as the one for wave map equations, such as the $O(3)$-$\sigma$ model in dimension two.

More precisely, we consider the solution for the Poiseuille laminar flow via a tube for \eqref{fels} (or \eqref{wels} with $k_1=k_2$). To describe the Poiseuille flow, let $(r,\theta, z)$ denote the cylindrical
coordinates of $\mathbb R^3$.
We consider a special form of solutions, so called Poiseuille type solutions. For $(x,y)\in \mathbb R^2$ and $r=\sqrt{x^2+y^2}$, consider  
$$
\u(x,y,t)=(0,0,v(r,t)),\quad P(x,y,t)=P(r,t)
$$
$$
\d(x,y,t)=\big(\sin\phi(r,t) \cos k\theta,\sin\phi(r,t)\sin k\theta, \cos\phi(r,t) \big).
$$
Then the system \eqref{fels} becomes
\begin{equation}\label{simeqn}
\begin{cases}
\displaystyle v_t=\frac{1}{r}\Big(rv_r+r\phi_t\Big)_r,\\
\\
 \displaystyle \phi_{tt}+2\phi_t=\frac{1}{r}\big(r\phi_r\big)_r-k^2\frac{\sin(2\phi)}{2r^2}-v_r.
\end{cases}
\end{equation}
We consider the following initial and boundary values
\beq\label{initial} 
v(r,0)=v_0,\quad \phi(r,0)=\phi_0,\quad \partial_t\phi(r,0)=\phi_1,
\eeq
\beq\label{bdycon}
v(0,t)=0,\quad \phi(0,t)=0,\quad \phi(\infty, t)=\pi.
\eeq

The main result of this paper is the following theorem.

\btm\label{mthmain} There exist initial values $v_0$, $\phi_0$ and $\phi_1$, such that the local smooth solutions to the initial-boundary value problem \eqref{simeqn}-\eqref{bdycon} blow up at finite time.
\etm

%

\begin{remark}
The detail of initial values and description of singularity will be given in Theorem \ref{mthsing}.
\end{remark}

One notices that the second equation of \eqref{simeqn} includes a damped
$O(3)$-$\sigma$ model. If we omit the coupling on $v$ and the damping term $2\phi_t$, it becomes the axisymmetric $O(3)$-$\sigma$ model
\beq\label{o3si}
\phi_{tt}=\frac{1}{r}\big(r\phi_r\big)_r-k^2\frac{\sin(2\phi)}{2r^2},\qquad k\geq 1,
\eeq
which can be derived from the wave map equation under a symmetry in dimension two. See \cite{[40]} for the model and existence of finite energy ($H^1$ type) solution.

Let's review some fundamental works for the $O(3)$-$\sigma$ model. For smooth initial data $(\phi_0,\phi_1)$, existence of global regular solutions to \eqref{o3si} has been established in the pioneering works by Shatah-Tahvildar-Zadeh \cite{[41]} by assuming the initial energy $E(\phi)(0)$ is sufficiently small (see also \cite{ssb98} Theorem 8.1), where 
\[
\tilde E(\phi)(t)=\pi\int_{\mathbb R^+}\left(
(\phi_t)^2+(\phi_r)^2+\frac{k^2}{r^2}\sin^2 \phi\right)\, rdr.
\]


In \cite{St}, Struwe showed that the singularities of equation \eqref{o3si} can only be developed by concentrating energy at the tip of a light cone by bubbling off at least one non-trivial harmonic map with finite energy. 

When $\phi(0,t)=0,\ \phi(\infty, t)=\pi,$
there are several approaches establishing examples of singularity formation of \eqref{o3si} with $\tilde{E}(\phi)(0)=4k \pi+\varepsilon$: Rodnianski and Sterbenz \cite{RS} when $k\geq 4$, Krieger-Schlag-Tataru \cite{CST} when $k=1$, and \cite{RR} for all other cases. 

For solutions satisfying $\phi(0,t)=0,\ \phi(\infty, t)=0,$
 the global existence and scattering have been proved by a sequence of seminal works: C\^{o}te-Kenig-Merle \cite{CKM},
C\^{o}te-Kenig-Lawrie-Schlag \cite{CKL},  Jendrej-Lawrie \cite{JL}, when $$\tilde{E}(\phi)(0)\leq 8k \pi.$$
There are many other fundamental works on the existence and behavior of solutions for \eqref{o3si}, such as \cite{SU}.

Our main analysis builds on the framework by \cite{RS} for \eqref{o3si} with $k\geq4$. One major difficulty in studying the singularity formation of \eqref{simeqn} is that there is lack of an $H^1$ estimate on velocity $v$ using the standard energy method (Lemma \ref{englemma1}). Inspired by the idea for the one dimensional quasilinear model in \cite{CHL20}, we introduce a new variable $h(r,t)$ defined in \eqref{defh} and satisfying \eqref{eqnv}. Using the variables $(\phi,h)$, we can show an improved regularity in Lemma \ref{englemma2}, which provides the desired weighted $H^1$ estimate on $v$. Here the improved regularity is not simply obtained by the heat equation in \eqref{simeqn}. Instead it is proved using the different scales in heat ($\partial_t\approx \partial_{xx}$) and wave equations ($\partial_{tt}\approx \partial_{xx}$).
The strong coupling of the heat and wave map equations also cause major difficulties in our analysis.


The remaining of the paper is organized as follows. In Section 2, we prove several a priori energy estimates for the system \eqref{simeqn}. In Section 3, we review the notations of the system and give the required initial data. A detailed main result on singularity formation is given at the end of this section. In Section 4, we derive the first order equation of the trajectory $\lambda(t)$, and prove the existence and orbital stability of $\lambda(t)$. In Section 5 and Section 6, we derive the second order equation of $\lambda(t)$ and prove the key technique estimates by bootstrapping arguments. In Section 7, we show the existence singularity at finite time. In Appendix A, we derive the simplified system \eqref{simeqn} from general Ericksen-Leslie system via the Poiseuille flow with special coefficients. In Appendix B, we show a stronger decay estimate outside a sufficiently large cone.


\section{Energy estimates}
\setcounter{equation}{0}

In this section, we provide several energy estimates to the system \eqref{simeqn} with initial and boundary values \eqref{initial} and \eqref{bdycon}.  
Without confusion, in this paper, we always use $C$ to denote the positive constant, which may take different values. 

\begin{lemma}\label{englemma1}
Any smooth solution to the system \eqref{simeqn} with initial and boundary values \eqref{initial} and \eqref{bdycon} satisfies the following energy inequality
\beq\label{enginq1}
\begin{split}
&\frac{d}{dt}\frac12\int_{0}^{\infty}\left(|v|^2+|\phi_t|^2+|\phi_r|^2+k^2\frac{\sin^2\phi}{r^2}\right)r\,dr\\
=&-\int_{0}^{\infty}(2-c_0^2)|\phi_t|^2r\,dr-\int_{0}^{\infty}\left(1-\frac{1}{c_0^2}\right)|v_r|^2r\,dr
-\int_{0}^{\infty}\left(c_0\phi_t-\frac1{c_0}v_r\right)^2r\,dr,
\end{split}
\eeq
where $c_0$ is a constant with $1<c_0^2<2$.
\end{lemma}

\begin{remark}
For later reference, we denote the energy decay by
\beq\label{Energy1}
E(t):=E[\phi,v](t)=\pi\displaystyle\int_{\mathbb{R}^+}[(\phi_t)^2+(\phi_r)^2+\frac{k^2}{r^2}\sin^2(\phi)+v^2]rdr\leq E[\phi,v](0)=:E(0).\eeq
\end{remark}

\pf Multiplying the first equation of \eqref{simeqn} by $vr$ then integrating it over $(0,\infty)$, we obtain
\beq\notag
\begin{split}
\frac{d}{dt}\frac12\int_{0}^{\infty}|v|^2r\,dr
=-\int_{0}^{\infty}|v_r|^2r\,dr-\int_{0}^{\infty}\phi_tv_rr\,dr.
\end{split}
\eeq
Multiplying the second equation of \eqref{simeqn} by $\phi_tr$ and integrating over $(0,\infty)$, we obtain
\beq\notag
\begin{split}
\frac{d}{dt}\frac12\int_{0}^{\infty}\left(|\phi_t|^2+|\phi_r|^2+k^2\frac{\sin^2\phi}{r^2}\right)r\,dr
=-2\int_{0}^{\infty}|\phi_t|^2r\,dr-\int_{0}^{\infty}\phi_tv_rr\,dr.
\end{split}
\eeq
Adding above two equations together, we conclude
\beq\label{engpf1}
\frac{d}{dt}\frac12\int_{0}^{\infty}\left(|v|^2+|\phi_t|^2+|\phi_r|^2+k^2\frac{\sin^2\phi}{r^2}\right)r\,dr
=-\int_{0}^{\infty}\left(2|\phi_t|^2+|v_r|^2\right)r\,dr-2\int_{0}^{\infty}\phi_tv_rr\,dr.
\eeq
For the right side, we may take a positive constant $c_0$ such that
\beq\notag
1<c_0^2<2
\eeq
and then it holds
\beq\notag
2|\phi_t|^2+|v_r|^2-2\phi_tw_r
=(2-c_0^2)|\phi_t|^2+\left(1-\frac{1}{c_0^2}\right)|v_r|^2
+\left(c_0\phi_t-\frac1{c_0}v_r\right)^2.
\eeq
Plugging it into \eqref{engpf1}, it holds
\beq\notag
\begin{split}
&\frac{d}{dt}\frac12\int_{0}^{\infty}\left(|v|^2+|\phi_t|^2+|\phi_r|^2+k^2\frac{\sin^2\phi}{r^2}\right)r\,dr\\
=&-\int_{0}^{\infty}(2-c_0^2)|\phi_t|^2r\,dr-\int_{0}^{\infty}\left(1-\frac{1}{c_0^2}\right)|v_r|^2r\,dr
-\int_{0}^{\infty}\left(c_0\phi_t-\frac1{c_0}v_r\right)^2r\,dr,
\end{split}
\eeq
which completes the proof of \eqref{enginq1}.
\endpf

\medskip
We also need some higher order estimates of velocity, which  will play a crucial part in our construction of singularity formations. The idea is inspired by arguments in \cite{CHL20} for one dimensional problem. The main observation is that the scales of $t$ and $r$ are different in heat and wave equations. More precisely, we introduce a new variable
\beq\label{defh}
h(r,t)=\frac1r\int_0^rv(R,t)\,R\, dR.
\eeq
To obtain the equation of $h$, multiplying the first equation of \eqref{simeqn} by $r$, then integrating it over $(0,r)$ and using the fact $(rv_r+r\phi_t)\big|_{r=0}=0$, we have
$$
\left(\int_0^rv(R,t)\,R\, dR\right)_t=rv_r+r\phi_t.
$$
Direct computation implies that
$$
(rh)_r=rv,\quad v_r=\left(\frac{(rh)_r}{r}\right)_r
=h_{rr}+\frac{h_r}{r}-\frac{h}{r^2}.
$$
The equation of $h$ is as follows
\beq\label{eqnv}
h_t=v_r+\phi_t=\frac1r(rh_r)_r-\frac{h}{r^2}+\phi_t.
\eeq
By the definition, the boundary condition is $h(0,t)=0$.
The second equation in \eqref{simeqn} can be written as follows
 \beq\label{eqnva}
  \phi_{tt}+\phi_t=\frac{1}{r}\big(r\phi_r\big)_r-k^2\frac{\sin(2\phi)}{2r^2}-h_t.
 \eeq

\begin{lemma}\label{englemma2}
Any smooth solution $(h, \phi)$ to the system \eqref{eqnv} and \eqref{eqnva} with initial and boundary values \eqref{initial} and \eqref{bdycon} satisfies the following energy estimates
  \beq\label{2.17}
 \begin{split}
\frac{d}{dt}\int_0^\infty\left(|h_t|^2+|h_r|^2+\frac{h^2}{r^2}\right)r\,dr
+\int_0^\infty\left(|h_{tr}|^2+\frac{|h_t|^2}{r^2}+|h_t|^2\right)r\,dr
\leq \frac{C}{\pi}E(0).
 \end{split}
 \eeq
\end{lemma}

\pf
Multiplying \eqref{eqnv} by $h_tr$ and integrating it over $(0,\infty)$ with respect to $r$, we obtain
 \beq\notag
 \frac{d}{dt}\frac{1}{2}\int_0^{\infty}\left(|h_r|^2+\frac{h^2}{r^2}\right)r\,dr+\int_0^{\infty}|h_t|^2r\,dr
 =\int_0^{\infty}h_t\phi_tr\,dr.
 \eeq
Applying the Young inequality, we have
  \beq\label{engpf2}
 \frac{d}{dt}\int_0^{\infty}\left(|h_r|^2+\frac{h^2}{r^2}\right)r\,dr+\int_0^{\infty}|h_t|^2r\,dr
 \leq C\int_0^{\infty}|\phi_t|^2r\,dr\leq \frac{C}{\pi}E(0).
 \eeq
 To obtain higher order estimates, we may take $t$-derivative for both sides of equation of $h$
 \beq\label{eqnvt}
 \begin{split}
 h_{tt}=&\frac1r(rh_{tr})_r-\frac{h_t}{r^2}+\phi_{tt}\\
 =&\frac1r(rh_{tr})_r-\frac{h_t}{r^2}-\phi_t+\frac1r(r\phi_r)_r-\frac{k^2\sin(2\phi)}{2r^2}-h_t,
 \end{split}
 \eeq
where we have used the equation \eqref{eqnva} in the last step.
Multiplying \eqref{eqnvt} by $rh_t $, then integrating by parts and using the fact $h_t(0,t)=0$, we obtain
 \beq\notag
 \begin{split}
\frac{d}{dt}\frac12\int_0^\infty|h_t|^2r\,dr
=-\int_0^\infty\left[|h_{tr}|^2+\frac{|h_t|^2}{r^2}+|h_t|^2
+\left(\phi_t-\frac1r(r\phi_r)_r+\frac{k^2\sin(2\phi)}{2r^2}\right)h_t\right]r\,dr.
 \end{split}
 \eeq
 For the last term, we estimate it as follows
 \beq\notag
 \begin{split}
\left| \int_0^\infty\left(\phi_t+\frac{k^2\sin(2\phi)}{2r^2}\right)h_tr\,dr\right|
\leq \frac{1}{2}\int_0^\infty\left(|h_t|^2+\frac{|h_t|^2}{r^2}\right)r\,dr
+C\int_0^\infty\left(|\phi_t|^2+\frac{\sin^2\phi}{r^2}\right)r\,dr
\end{split}
 \eeq
 \beq\notag
 \left|\int_0^\infty\frac1r(r\phi_r)_rh_tr\,dr\right|
 =\left|\int_0^\infty \phi_rh_{tr}r\,dr\right|
 \leq \frac12\int_0^\infty|h_{tr}|^2r\,dr+C\int|\phi_r|^2r\,dr.
 \eeq
 Therefore
  \beq\label{engpf3}
 \begin{split}
\frac{d}{dt}\int_0^\infty|h_t|^2r\,dr
+\int_0^\infty\left(|h_{tr}|^2+\frac{|h_t|^2}{r^2}+|h_t|^2\right)r\,dr
\leq C\int_0^\infty\left(|\phi_t|^2+|\phi_r|^2+\frac{\sin^2\phi}{r^2}\right)r\,dr,
 \end{split}
 \eeq
 which combining with \eqref{engpf2} completes the proof of \eqref{2.17}.
 \endpf
 
 \medskip
By the estimate \eqref{Energy1} and the equation of $h_t$ in \eqref{eqnv}, one can find $L^2$ bound on $v_r$. But it is much more convenient to directly use $(h,\phi)$ as dependent variables than using $(v,\phi)$. One clear advantage is there is an extra estimate on $h_{tr}$ which will help us get the following higher order energy estimates with weights.
 
\begin{lemma} Let $(h,\phi)$ be any smooth solution to the system \eqref{eqnv} and \eqref{eqnva} with initial and boundary values \eqref{initial} and \eqref{bdycon}. Assume that the initial data satisfy 
\beq\label{init2}
\int_0^\infty\left[|h_t|^2+|\phi_t|^2+|\phi_r|^2+|h_{tr}|^2\right](r,0)\,r^{2+\delta}\,dr\leq \tilde C
\eeq
for some positive constant $\tilde C$, then it holds
 \beq\label{2.32}
 \int_0^\infty\left[|h_t|^2+|\phi_t|^2+|\phi_r|^2\right](r,t)\,r^{2+\delta}\,dr\leq C,
 \eeq
 and furthermore
 \beq\label{2.33}
 \int_0^t\int_0^\infty|h_{tr}|^2r^{2+\delta}\,drds\leq C,
 \eeq
for any $t\in[0,1]$. 
\end{lemma}

\begin{remark} From now on, we restrict our consideration in $t\in[0,1]$. Later, we will choose the initial data satisfying \eqref{init2} such that the singularity forms before $t=1$, to make the proof self-consistent.
\end{remark}

\pf
Multiplying \eqref{eqnvt} by $r^{2+\delta}h_t $ for any $0<\delta<1$ and integrating over $(0,\infty)$ with respect to $r$, we obtain
 \beq\label{engpf4}
\begin{split}
&\frac{d}{dt}\frac12\int_0^\infty|h_t|^2r^{2+\delta}\,dr\\
=&\int_0^\infty(rh_{tr})_rh_tr^{1+\delta}\,dr-\int_0^\infty\left[|h_t|^2r^{\delta}+|h_t|^2r^{2+\delta}\right]\,dr\\
&-\int_0^\infty
\left(\phi_t-\frac1r(r\phi_r)_r+\frac{k^2\sin(2\phi)}{2r^2}\right)h_tr^{2+\delta}\,dr.
 \end{split}
 \eeq
For the first term of right side, integrating by parts and using the fact $h_t(0,t)=0$, it holds
 \beq\notag
 \begin{split}
\int_0^\infty(rh_{tr})_rh_tr^{1+\delta}\,dr
=-\int_0^\infty |h_{tr}|^2r^{2+\delta}\,dr-(1+\delta)\int_0^\infty h_{tr}h_tr^{1+\delta}\,dr.
 \end{split}
 \eeq
 For the extra term, we may estimate it as follows 
  \beq\notag
 \begin{split}
&\left|\int h_{tr}h_tr^{1+\delta}\,dr\right|\\
=&\left|\int_0^1h_{tr}h_tr^{1+\delta}\,dr+\int_1^{\infty} h_{tr}h_tr^{1+\delta}\,dr\right|\\
\leq&\int_0^1|h_{tr}||h_t|r\,dr
+\xi\int_1^{\infty}|h_{tr}|^2r^{2+\delta}\,dr
+C\int_1^{\infty}|h_t|^2r^{\delta}\,dr\\
\leq& \xi\int_0^\infty|h_{tr}|^2r^{2+\delta}\,dr
+C\int_0^\infty\left[|h_{tr}|^2+|h_t|^2\right]r\,dr,
 \end{split}
 \eeq
 for some constant $\xi>0$, which will be determined later. 
Similarly, we can estimate the third term of right hand side of \eqref{engpf4} as follows
   \beq\notag
 \begin{split}
 \left|\int_0^\infty\phi_th_tr^{2+\delta}\,dr\right|
 \leq \frac12\int_0^\infty\left[|\phi_t|^2+|h_t|^2\right]r^{2+\delta}\,dr,
 \end{split}
 \eeq
   \beq\notag
 \begin{split}
&\left|\int_0^\infty (r\phi_r)_rh_tr^{1+\delta}\,dr\right|\\
=& \left|\int_0^\infty\left[\phi_rh_{tr}r^{2+\delta}+(1+\delta)\phi_rh_tr^{1+\delta}\right]\,dr\right| \\
\leq&\xi\int_0^\infty|h_{tr}|^2r^{2+\delta}\,dr
+\xi\int_0^\infty|h_{t}|^2r^{\delta}\,dr
+C\int_0^\infty|\phi_r|^2r^{2+\delta}\,dr,
 \end{split}
 \eeq
 and
   \beq\notag
 \begin{split}
 &\left|\int_0^\infty\frac{\sin(2\phi)}{2r^2}h_tr^{2+\delta}\,dr\right|\\
 =&\left|\int_0^1\sin(2\phi)h_tr^{\delta}\,dr+\int_1^\infty\sin(2\phi)h_tr^{\delta}\,dr\right|\\
 \leq& C\int_0^\infty\left[\frac{\sin^2\phi}{r^2}+|h_t|^2\right]r\,dr
 +C\int_0^\infty|h_t|^2r^{2+\delta}\,dr.
 \end{split}
 \eeq
 Combining all above estimates together with \eqref{engpf4}, we have
 \beq\label{2.24}
 \begin{split}
&\frac{d}{dt}\frac12\int_0^\infty|h_t|^2r^{2+\delta}\,dr
+(1-2\xi)\int_0^\infty|h_{tr}|^2r^{2+\delta}\,dr+(1-\xi)\int_0^\infty|h_t|^2r^{\delta}\,dr\\
\leq &C+C\int_0^\infty|h_{tr}|^2r\,dr+C\int_0^\infty\left[|\phi_t|^2+|h_t|^2\right]r^{2+\delta}\,dr.
 \end{split}
 \eeq
 To estimate the last term related to $\phi$, we need to use the equation \eqref{eqnva}. Multiplying the equation \eqref{eqnva} by $\phi_tr^{2+\delta}$ and integrating over $(0,\infty)$ with respect to $r$, we have
 \beq\label{2.25}
 \begin{split}
 &\frac{d}{dt}\frac12\int_0^\infty|\phi_t|^2r^{2+\delta}\,dr+\int_0^\infty|\phi_t|^2r^{2+\delta}\,dr\\
=&\int_0^\infty(r\phi_{r})_r\phi_tr^{1+\delta}\,dr-\int_0^\infty
\frac{k^2\sin(2\phi)}{2}\phi_tr^{\delta}\,dr-\int_0^\infty h_t\phi_tr^{2+\delta}\,dr.
 \end{split}
 \eeq
 Similar as the above arguments, we can estimate the first term of right hand side as follows
 \beq\notag
 \begin{split}
 \int_0^\infty(r\phi_{r})_r\phi_tr^{1+\delta}\,dr
 =-\frac{d}{dt}\frac12\int_0^\infty|\phi_r|^2r^{2+\delta}\,dr
 -(1+\delta)\int_0^\infty\phi_r\phi_tr^{1+\delta}\,dr.
 \end{split}
 \eeq
 The last term can be estimated as follows
 \beq\notag
 \begin{split}
 &\left|\int_0^\infty\phi_r\phi_tr^{1+\delta}\,dr\right|\\
 =&\left|\int_0^1\phi_r\phi_tr^{1+\delta}\,dr
 +\int_1^\infty\phi_r\phi_tr^{1+\delta}\,dr\right|\\
 \leq &C\int_0^\infty\left[|\phi_r|^2+|\phi_t|^2\right]r\,dr+C\int_0^\infty|\phi_t|^2r^{2+\delta}\,dr.
 \end{split}
 \eeq
 For the other terms on right hand side of \eqref{2.25}, we have
   \beq\notag
\left|\int_0^\infty\frac{\sin(2\phi)}{2}\phi_tr^{\delta}\,dr\right|
 \leq C\int_0^\infty\left[\frac{\sin^2\phi}{r^2}+|\phi_t|^2\right]r\,dr
 +C\int_0^\infty|\phi_t|^2r^{2+\delta}\,dr,
 \eeq
 and
   \beq\notag
 \begin{split}
 \left|\int_0^\infty\phi_th_tr^{2+\delta}\,dr\right|
 \leq \frac12\int_0^\infty\left[|\phi_t|^2+|h_t|^2\right]r^{2+\delta}\,dr.
 \end{split}
 \eeq
 Combining all above estimates together with \eqref{2.25}, we have
  \beq\label{engpf5}
 \frac{d}{dt}\frac12\int_0^\infty\left[|\phi_t|^2+|\phi_r|^2\right]r^{2+\delta}\,dr
\leq C+C\int_0^\infty\left[|\phi_t|^2+|h_t|^2\right]r^{2+\delta}\,dr.
 \eeq
 Therefore, combining estimates \eqref{2.24} and \eqref{engpf5}, choosing $\xi$ small enough and using \eqref{Energy1} and \eqref{2.17}, we finally obtain
 \beq\notag
 \begin{split}
\frac{d}{dt}\int_0^\infty\left[|h_t|^2+|\phi_t|^2+|\phi_r|^2\right]r^{2+\delta}\,dr
\leq C+C\int_0^\infty\left[|h_t|^2+|\phi_t|^2+|\phi_r|^2\right]r^{2+\delta}\,dr,
 \end{split}
 \eeq
then standard Gronwall arguments imply that
 \beq\notag
 \int_0^\infty\left[|h_t|^2+|\phi_t|^2+|\phi_r|^2\right]r^{2+\delta}\,dr\leq C,
 \eeq
for any $t\in[0,1]$, where we use the initial assumption \eqref{init2}. Then using \eqref{2.24}, \eqref{2.32} and \eqref{init2}, we can easily prove that
 \beq\notag
 \int_0^t\int|h_{tr}|^2r^{2+\delta}\,drds\leq C,
 \eeq
 for any $t\in[0,1]$, which completes the proof of lemma.
 \endpf

\section{Basic setup and main results of singularity formation}
\setcounter{equation}{0}
In the rest of this paper, we will establish a singularity formation for system \eqref{simeqn} which is also equivalent to the following system
\beq\label{vphi}
\begin{split}
h_t&=\frac1r(rh_r)_r-\frac{h}{r^2}+\phi_t,
\\
\phi_{tt} +\phi_t&=\frac{1}{r}(r\phi_r)_r-k^2\frac{\sin (2\phi )}{2r^2}-h_t,
\end{split}
\eeq
where $h$ is defined as in \eqref{defh}. Our construction is inspired by the idea of \cite{RS}. We will first review the basic setup as in \cite {RS} for the readers' convenience.


First notice that we can rewrite the energy functional as follows
\beq\label{sec3pf1}
\begin{split}
E[\phi,v]&=\displaystyle\pi\int_{\mathbb{R}^+}[(\partial_t\phi)^2+(\partial_r\phi)^2+\frac{k^2}{r^2}\sin^2(\phi)+v^2]rdr\\
&=\displaystyle\pi\int_{\mathbb{R}^+}[(\partial_t\phi)^2+(\partial_r\phi-\frac{k}{r}\sin(\phi))^2+v^2]rdr+2\pi\displaystyle\int_{\mathbb{R}^+}[k\sin(\phi)\partial_r\phi]dr\\
&=\displaystyle\pi\int_{\mathbb{R}^+}[(\partial_t\phi)^2+(\partial_r\phi-\frac{k}{r}\sin(\phi))^2+v^2]rdr+4k\pi.
\end{split}
\eeq
Let $I(r)$ and $v\equiv 0$ be the stationary solution to the energy functional \eqref{sec3pf1} with boundary data \eqref{bdycon}. Then $I(r)$ satisfies 
\beq\label{I_def}
r\partial_r I =k\sin(I):=J.
\eeq
It is not hard to solve
\beq\label{I-def2}
I=2\arctan (r^k).
\eeq
For a scaling trajectory $\lambda(t)$ which will be defined later,
we denote that
\beq\label{I_def2}
I_\lambda(r,t)=I(\lambda(t) r)=2\arctan (\lambda r)^k,\quad 
J_\lambda(r,t)=k\sin(I_\lambda)=r\partial_r I_\lambda=r\partial_rI(\lambda r).
\eeq
For convenience, we denote 
\[
\dot{I}_\lambda=\frac{\partial}{\partial t}(I_\lambda),\qquad
\ddot{I}_\lambda=\frac{\partial^2}{\partial{t^2}}(I_\lambda),
\]
so does $\dot{J}_\lambda$ and $\ddot{J}_\lambda$.
It is easy to check that 
\beq\label{jte}
\dot{I}_\lambda=\frac{\dot{\lambda}}{\lambda}J_\lambda,\quad
\dot{J}_\lambda=\frac{\dot{\lambda}}{\lambda}(r\partial_rJ)_\lambda.
\eeq

As  in \cite{RS}, the linearized Hamiltonian $H_\lambda$ is given as follows
\beq\label{hleq}
H_\lambda =-\partial _{rr}-\frac{1}{r}\partial_r+\frac{k^2}{r^2}\cos(2I_\lambda)=A_\lambda^*A_\lambda,
\eeq
where
\beq\notag
A_\lambda=-\partial_r+\frac{k}{r}\cos(I_\lambda),\quad 
A_\lambda^*=\partial_r+\frac1r+\frac{k}{r}\cos(I_\lambda).
\eeq
It is not hard to see that $J_\lambda$ solves the first order linearized Bogomol'nyi equation
 \beq\label{aljl}
 A_\lambda J_\lambda =0.
 \eeq
The conjugate operator of $H_\lambda$ is defined as follows
\beq\label{def_Htil}
\tilde{H}_\lambda=A_\lambda A_\lambda^*
=-\partial _{rr}-\frac{1}{r}\partial_r+V_\lambda(r),
\eeq 
where 
$$
V_\lambda(r)=\frac{k^2+1}{r^2}+\frac{2k}{r^2}\cos(I_\lambda).
$$
The following relation holds for $H_\lambda$ and $\tilde H_\lambda$
\beq
A_\lambda H_\lambda=\tilde H_\lambda A_\lambda.
\eeq
The following positive, space-repulsive and time-repulsive properties hold for $V_\lambda$ (cf. (51)-(53) in \cite{RS})
\beq\label{postives}
V_\lambda\geq \frac{(k-1)^2}{r^2},\quad -\partial_rV_\lambda\geq 2\frac{(k-1)^2}{r^3},\quad -\partial_tV_\lambda=\frac{2k^2}{r^2}\frac{\dot{\lambda}}{\lambda}\sin^2(I_\lambda).
\eeq

Denote 
 $$u:=\phi(t,r)-I_\lambda$$ 
which is equivalent to the following decomposition of $\phi$
\beq\label{phiiu}
\phi(t,r)=I_\lambda +u.
\eeq
By \eqref{vphi} and direct computation, we obtain the equation of $u$ as follows


\beq\label{3.39}
u_{tt}+u_t+H_\lambda u=-\ddot{I}_\lambda-\dot{I}_\lambda+\mathcal{N}(u)-h_t,
\eeq
where
\beq\notag
\mathcal{N}(u)=\frac{k^2\sin(2I_\lambda)}{2r^2}(1-\cos(2u))+\frac{k^2\cos(2I_\lambda)}{r^2}(u-\frac{1}{2}\sin(2u)).
\eeq

Furthermore, we may decompose $u$ as follows (cf. \cite{RS})
\beq\label{uw0w}
u=w_0+w
\eeq
where
\beq\label{defw0}
w_0(t,r)=\frac{\dot{\lambda}^2}{\lambda^4}(t)(aJ_\lambda(r)+b(r^2J)_\lambda),
\eeq 
with
$a=-\frac{1}{4}\langle J,r^2J\rangle\cdot \|J\|^{-2}_{L^2(rdr)}$, $b=\frac{1}{4}.$ 
By the definition of $w_0$, one may claim that 
\beq\label{w0AH}
\langle w_0,J_\lambda\rangle =0,\quad A_\lambda H_\lambda(w_0)=\tilde H_\lambda A_\lambda(w_0)=-A_\lambda(\ddot{I}_\lambda+2\dot{I}_\lambda).
\eeq
In fact, direct computation implies the orthogonality. To show the second equation, first notice that
\beq\label{AId1}
-A_\lambda\dot{I}_\lambda=-\frac{\dot{\lambda}}{\lambda}A_\lambda(r\partial_rI_\lambda)=-\frac{\dot{\lambda}}{\lambda}A_\lambda J_\lambda=0,
\eeq
\beq\notag
A_\lambda\ddot{I}_\lambda=\frac{\dot{\lambda}}{\lambda}A_\lambda\partial_t(r\partial_r I_\lambda)=\frac{\dot{\lambda^2}}{\lambda^2}A_\lambda(r\partial_rJ_\lambda)=\frac{\dot{\lambda}^2}{\lambda^2}A_\lambda(aJ+r\partial_rJ)_\lambda.
\eeq
Combining these with the definition of $w_0$, 
one can verify the second equation in \eqref{w0AH}. 

\medskip
Following the notation in \cite{RS}, for any pair of nonnegtive integers $m,n$, $F^{m,n}$ denotes any $C^{\infty}(\R^+)$ function with the bounds
$$
|(r\partial_r)^i F^{m,n}|\leq C\frac{r^m}{(1+r)^{m+n}},
$$
for any nonnegative integer $i$. For simplicity, denote $F^{n}=F^{0,n}$. Also, denote $F^{m,n}_\lambda(r)=F^{m,n}(\lambda r)$ as the rescaling of $F^{m,n}$. These functions are very useful to record the order of $r$ in the estimates and may be different from line to line. By the definition, it is not hard to see 
$$
r^j F^{m,n}_\lambda=\lambda^{-1}F^{m+j,n-j}_\lambda,\quad
\partial_tF^{m,n}_\lambda=\dot{\lambda}\lambda^{-1}F^{m,n}_\lambda,\quad 
\partial_r F^{m,n}_\lambda=\lambda F^{m-1,n+1}_\lambda
$$
where the integer $j$ satisfies $j\geq-m$ in first identity and $m\geq 1$ in third identity. 

Applying the notation with restriction $k\geq 4$, we obtain
\beq\label{58a1}
w_0=\dot{\lambda}^2\lambda^{-4}(F^{4,4}_{\lambda}+F^{6,2}_{\lambda}),
\eeq
\beq\label{58a2}
\partial_t w_0=\ddot{\lambda}\dot{\lambda}\lambda^{-4}(F^{4,4}_{\lambda}+F^{6,2}_{\lambda})+\dot{\lambda}^3\lambda^{-5}(F^{4,4}_{\lambda}+F^{6,2}_{\lambda}),
\eeq
\beq\label{58a3}
A_\lambda w_0=\dot{\lambda}^2\lambda^{-3}(F^{3,5}_{\lambda}+F^{5,3}_{\lambda}),
\eeq
\beq\label{58b}
|J_\lambda|\leq C F^{4,4}_\lambda,\quad
|\partial_t J_\lambda|\leq C \dot{\lambda}\lambda^{-1}F^{4,4}_{\lambda},\quad 
|\partial^2_t J_\lambda|\leq C (|\ddot{\lambda}|\lambda^{-1}+\dot{\lambda}^2\lambda^{-2})F^{4,4}_{\lambda},
\eeq
\beq\label{58c}
A_\lambda(F^{m,n}_\lambda)=\lambda F^{m-1,n+1}_\lambda,\quad
\partial_t A_{\lambda}=\dot{\lambda}F^{7,9}_\lambda,\quad
|\partial^2_t A_\lambda|\leq C(\ddot{\lambda}+\dot{\lambda}^2\lambda^{-1})F^{7,9}_{\lambda},
\eeq
where we need $m\geq 1$ in the first identity of \eqref{58c}.

\subsection{Main results}
For simplicity, we assume the initial data
\beq\label{vint0}
v_0\equiv 0
\eeq
which implies that $h_0(r):=h(r,0)\equiv 0$ for the system \eqref{vphi}. We also assume 

\beq\label{phi00}
\phi_0(r)=u_0(r)+I(\lambda_0r),
\quad 
\phi_1(r)=\la_1J(\lambda_0r)+g_0(r)
\eeq
where 
\beq\label{lambda1}
\lambda_0=\v^{-4},\quad \lambda_1=\frac{\dot\la}{\la}(0)=C_0^2\frac{\v^5}{\pi}\|J(\lambda_0r)\|_{L^2(rdr)}^{-2}
\eeq for some small constant $\v>0$, $g_0(r)$ is a smooth given function and 
\beq\label{u0per}
\int_{\R^+}u_0\cdot J(\lambda_0r)\,rdr=0.
\eeq
 For another sufficiently small positive constant $c_0$ with $c_0^2\geq \v$, we assume that
 \beq\label{2930}
 \sum_{i=0}^{1}\int_{\R^+}(1+r^2)^{1-i}
 \left(
 (\partial^{i+1}_ru_0)^2+\frac{(\partial^i_ru_0)^2}{r^2}
 +(\partial^{i}_rg_0)^2+\frac{g_0^2}{r^2}
 \right)\,rdr\leq c_0^2\v^2.
 \eeq
Notice that it is not hard to check that \eqref{init2} is always satisfied, hence \eqref{2.32} and \eqref{2.33} hold as long as the solution is smooth.

Now we are ready to state our main result on the singularity formation.
 
 
\btm\label{mthsing} Suppose that $k\geq 4$.
There exists a continuous time dependent parameter $\lambda(t)$ with $\lambda(0)=\lambda_0=\v^{-4}$, such that the regular solution $(h(r,t),\phi(r,t))$ to the system \eqref{vphi} with initial data \eqref{vint0}-\eqref{2930} and boundary data \eqref{bdycon} can be split into 
\beq\label{31}
\phi(r,t)=I(\lambda(t)r)+u(r,t).
\eeq
The reminder term $u(r,t)$ satisfies 
\beq\label{32}
\int_{\R^+}\left(u_t^2+u_r^2+\frac{u^2}{r^2}\right)\,rdr\leq  \v^2
\eeq 
for as long as the solution exists. There also exists a finite time $0<T^*\leq C \v^{\frac12}$ such that 
\beq
\lim\limits_{t\rightarrow T^*}\lambda(t)=\infty,
\eeq
and $\lambda(t)$ satisfies the following estimates
\beq\label{33}
\frac{M}{T^*-t}\leq \lambda(t)\leq c_0^{\frac14} \frac{\sqrt{|\ln(T^*-t)|}}{T^*-t}
\eeq
when $t$ is sufficiently close to $T^*$.
\etm


\section{Trajectory $\lambda(t)$ and orbital stability}
\label{sec4}
\setcounter{equation}{0}
This section will be devoted to study the existence and properties of the trajectory $\lambda(t)$. Let's first derive the equation of $\lambda(t)$ by assuming the orthogonality between $u$ and $J_\lambda$.

\begin{lemma}\label{lemma4.1} Let $u$ be defined as in \eqref{phiiu} and $\phi$ be sufficiently smooth. Suppose that the orthogonality 
\beq{\label{3.1}}
\langle u,J_\lambda \rangle=0,
\eeq
is valid for any time $t>0$. Then the trajectory $\lambda(t)$ satisfies the following equation
\beq\label{63}
\dot{\lambda}\left(2\langle I,J\rangle+\left\langle \phi\left(\frac{r}{\lambda}\right),r\partial_rJ\right\rangle\right)=-\langle\phi_t,J_\lambda\rangle \lambda^3.
\eeq
Reversely, if $\lambda(t)$ satisfies the equation \eqref{63} and the orthogonality is valid at initial time, then $u$ is orthogonal to $J_\lambda$ for any time $t>0$. 

\end{lemma}

\pf
Taking derivative with respect to time $t$ on equation \eqref{3.1}, we obtain
\beq\label{sec4pf1}
0=\frac{d}{dt}\langle u,J_\lambda\rangle =\langle u_t,J_\lambda\rangle+\langle u,\partial_tJ_\lambda\rangle\\
=\langle u_t,J_\lambda\rangle+\frac{\dot{\lambda}}{\lambda}\langle u,(r\partial_r J)_\lambda\rangle,
\eeq
where we use that \eqref{jte}.
Direct computation implies 
\[\begin{split}
\langle u_t,J_\lambda\rangle&=\langle \phi_t-(I_\lambda)_t,J_\lambda\rangle=\langle\phi_t,J_\lambda\rangle-\frac{\dot{\lambda}}{\lambda}\langle r\partial_rI_\lambda,J_\lambda\rangle=\langle \phi_t,J_\lambda\rangle-\frac{\dot{\lambda}}{\lambda^3}\langle J,J\rangle
\end{split}
\]
where we have used the definition of $J_\lambda$ and the fact 
$\langle J_\lambda,J_\lambda\rangle=\frac{1}{\lambda^2}\langle J,J\rangle$ in last step.  
Similarly, it holds
\[\begin{split}
&\frac{\dot{\lambda}}{\lambda}\langle u,(r\partial_rJ)_\lambda\rangle\\
=&\frac{\dot{\lambda}}{\lambda}\displaystyle\int_{\mathbb{R}^+}\phi(r\partial_rJ)_\lambda rdr-\frac{\dot{\lambda}}{\lambda}\displaystyle\int_{\mathbb{R}^+}I_\lambda(r\partial_rJ)_\lambda rdr\\
=&\frac{\dot{\lambda}}{\lambda^3}\displaystyle\int_{\mathbb{R}^+}\phi(r\partial_rJ)_\lambda \lambda rd\lambda r-\frac{\dot{\lambda}}{\lambda^3}\displaystyle\int_{\mathbb{R}^+}I(r\partial_rJ)rdr\\
=& \frac{\dot{\lambda}}{\lambda^3}\left\langle \phi\left(\frac{r}{\lambda}\right),r\partial_rJ\right\rangle+\frac{\dot{\lambda}}{\lambda^3}\big(2\langle I,J\rangle+\langle J,J\rangle\big)
\end{split}
\]
where we have used integration by parts in the last step.
Combining above equations with \eqref{sec4pf1} implies \eqref{63}. 
The section part of this lemma is trivial. This completes the proof of lemma.
\endpf


\medskip
The next result is on the existence and orbital stability of $\lambda(t)$ with suitable initial data.

%
%

\begin{proposition}\label{prop3.1} 
Let $(\phi,v)$ be any solution to the system \eqref{simeqn} with initial values \eqref{vint0}-\eqref{2930} and boundary values \eqref{bdycon}. As long as $(\phi,v)$ is smooth in time interval $(0,T]$, the equation \eqref{63} with initial data $\lambda(0)=\lambda_0>0$ has a continuous solution $0<\lambda(t)<\infty$, which also satisfies 
\beq\label{cond3}
\left|\frac{\dot{\lambda}}{\lambda^2}\right|\leq C\epsilon.
\eeq
Denote
$u=\phi(r,t)-I_\lambda$. 
Then both
\beq\label{cond1}
{E_0[u,v]=\int_{\mathbb{R}^+}\left[(\partial_t\phi)^2+(\partial_ru)^2+\frac{k^2}{r^2}u^2+v^2\right]rdr\leq C \epsilon^2
}
\eeq
and
\beq\label{cond2}
\langle u,J_\lambda \rangle=0
\eeq
hold in time interval $(0,T]$.

\end{proposition}

\bigskip
\pf  Now we prove the proposition in three steps. 

\paragraph{\bf Step 1.}

Denote $$\alpha(\lambda,t)=2\langle  I,J\rangle+
\left\langle\phi\left(\frac{r}{\lambda}\right),r\partial_rJ\right\rangle,\quad
\beta(\lambda,t)=-\langle\phi_t,J_\lambda\rangle\lambda^3.$$
The equation \eqref{63} becomes
\beq\label{sec4pf2}
\alpha(\lambda,t)\dot\lambda=\beta(\lambda,t).
\eeq
Since $\phi$ is smooth in time interval $(0,T)$, it is easy to see that
$ \alpha, \beta$ are $C^1$ functions of $\lambda$ and $t$.

\medskip
We first need to estimate $\alpha(\lambda_0,0)$. Notice that
\beq
\alpha(\lambda_0,0)+\langle J,J\rangle
=2\langle  I,J\rangle+\langle J,J\rangle+\left\langle\phi\left(\frac{r}{\lambda_0}\right),r\partial_rJ\right\rangle.
\eeq
Direct computation implies
$$\langle I,r\partial_rJ\rangle=\displaystyle\int_{\mathbb{R}^+}Ir\partial_rJrdr
=\displaystyle\int_{\mathbb{R}^+}Ir^2dJ=-2\langle I,J\rangle-\langle J,J\rangle,
$$
where we have used the facts: when $k\ge 3$
$$Ir^2J|_{r=0}=\left.2\arctan(r^k)rk\frac{2r^k}{1+r^{2k}}\right |_{r=0}=0,\quad
 Ir^2J|_{r=\infty}=\left.\pi k\frac{r^{k+2}}{1+r^{2k}}\right|_{r=\infty}=0.$$ 
Hence, it holds
\beq\label{sec4pf3}
\left.\begin{split}
 &|\alpha(\lambda_0,0)+\langle J,J\rangle|\\
=&|-\langle I,r\partial_rJ\rangle+\left\langle\phi\left(\frac{r}{\lambda_0},0\right),r\partial_rJ\right\rangle|\\
=&\left|\int_{\mathbb{R}^+}u\left(\frac{r}{\lambda_0},0\right)\partial_rJ\,r^2dr\right|\\
\leq& \left(\int_{\mathbb{R}^+}\left(\frac{u(\frac{r}{\lambda_0},0)}{r}\right)^2rdr\right)^{\frac{1}{2}}\left(\int_{\mathbb{R}^+} (r^2\partial_rJ)^2rdr\right)^{\frac{1}{2}}\\
=&\left(\int_{\mathbb{R}^+}\left(\frac{u_0}{r}\right)^2rdr\right)^{\frac{1}{2}}\left(\int_{\mathbb{R}^+} (\partial_rJ)^2r^5dr\right)^{\frac{1}{2}}
\le C\epsilon
\end{split}\right.
\eeq
where we use the assumption \eqref{2930} and the following fact when  $k\ge 3$
\beq\notag
\begin{split}
\int_{\mathbb{R}^+}(\partial_r J)^2r^5dr
=\int_{\mathbb{R}^+}\frac{4k^4r^{2k+3}(r^{2k}-1)^2}{(1+r^{2k})^4}dr
\le C+ 4k^4\displaystyle\int_{1}^\infty\frac{1}{r^{2k-3}}dr
<\infty.
\end{split}
\eeq
Choosing sufficiently small $\epsilon>0$ such that $C\epsilon$ is much smaller than 
\beq\label{C0}
\langle J,J\rangle =C_0,
\eeq
we obtain that $\alpha(\lambda_0,0)$ is bounded below with a positive lower bound. 

By the regularity of $\phi$ and the continuity of $\alpha(\lambda,t)$ and $\beta(\lambda,t)$ in $\lambda$, there exists a time  $T_1\in \mathbb{R}^+$ with $T_1\leq T$ such that the equation \eqref{63} has a continuous solution $\lambda(t)\geq 0$ in $[0,T_1]$, and for any $t\in [0,T_1]$
\beq\label{sec4pf4}
|\alpha(\lambda,t)+\langle J,J\rangle|\le 2C\epsilon.
\eeq
It is easy to see that the solution
$\lambda(t)$ of \eqref{63} also satisfies the orthogonality \eqref{cond2}
since it holds initially.

Since $\lambda(t)$ is continuous in $[0,T_1]$, we know $\lambda(t)<\infty$.  Suppose that $\lambda(t)\rightarrow 0^+$ at finite time. There will be a singularity of function $I(\lambda r)$ at $r=\infty$, which contradicts to the finite speed of propagation (this kind of singularity can only happen when $t=\infty$). Thus $\lambda>0$ for any finite time.

To show the solution $\lambda(t)$ exists for any time in $(0,T]$, we need to improve \eqref{sec4pf4} 
\beq\label{u-est-alpha}
|\alpha(\lambda,t)+\langle J,J\rangle|\le C\epsilon.
\eeq
for any time $t\in [0,T_1]$. By the similar arguments as \eqref{sec4pf3}, it is sufficient to show  \eqref{cond1} is valid for any time $t\in [0,T_1]$.

\paragraph{\bf Step 2.} 
In this step, we will show the estimate \eqref{cond1} is valid on $[0,T_1]$.

Recall, by the energy inequality \eqref{Energy1}, we have $E[\phi,v]\leq E(0)\leq 4k\pi+\epsilon^2$. Direct computation implies that
\beq\notag
\begin{split}
\epsilon^2\geq&E[\phi,v]-4k\pi=E[\phi,v]-E[I_\lambda,0]\\
=&\int_{\mathbb{R}^+}[(\partial_t\phi)^2+\left(\partial_r\phi-\frac{k}{r}\sin \phi\right)^2+v^2]rdr\\
=&\int_{\mathbb{R}^+}\left[(\partial_t\phi)^2+\left(\partial_ru-\frac{k}{r} \cos(I_\lambda)u\right)^2+v^2\right]rdr-\mathcal{R}_\lambda(u),
\end{split}
\eeq

where 
$$
\mathcal{R}_\lambda(u)=\int_{\mathbb{R}^+}\left(\partial_r\phi-\frac{k}{r}\sin\phi\right)^2rdr-\int_{\mathbb{R}^+}\left(\partial_ru-\frac{k}{r}\cos(I_\lambda)u\right)^2rdr.$$
Since $(\phi,v)$ is smooth, the energy $E_0[u, v]$ is finite.
By the estimate (162) in \cite{RS} Appendix B, it holds
$$\displaystyle\int_{\mathbb{R}^+}\left[(\partial_t\phi)^2+(\partial_ru)^2+\frac{k^2}{r^2}u^2\right]rdr
\le C\displaystyle\int_{\mathbb{R^+}}\left[(\partial_t\phi)^2+\left(\partial_ru-\frac{k}{r}\cos(I_\lambda)u\right)^2\right]rdr$$ 
which implies
\beq\label{est-E0R}
E_0[u,v]\leq C\epsilon^2+|\mathcal{R}_\lambda(u)|.
\eeq
By the definition of $I_\lambda$ and $u$, it holds
\[\partial_r \phi-\frac{k}{r}\sin(\phi)=\frac{k}{r}\sin I_\lambda+\partial_ru-\frac{k}{r}\sin(I_\lambda+u).
\]
Hence
$$
\mathcal{R}_\lambda(u)
=\int_{\mathbb{R}^+}
\left[2u_r+\frac{k}{r}\sin I_\lambda-\frac{k}{r}\sin(I_\lambda+u)-\frac{k}{r}\cos(I_\lambda)u\right]\left[\frac{k}{r}\sin I_\lambda-\frac{k}{r}\sin(I_\lambda+u)+\frac{k}{r}\cos(I_\lambda)u\right]\,rdr.$$
By the Taylor series expansion
\beq\notag
\sin(I_\lambda+u)=\sin(I_\lambda)+\cos(I_\lambda)u-\frac{1}{2}\sin(I_\lambda)u^2+\cdots
\eeq
it is not hard to conclude that 
\beq
|\mathcal{R}_\lambda(u)|
\leq C \int_{\mathbb{R}^+}\left(|u_r|+\frac{|u|}{r}\right)\left(\frac{u^2}{r}\right)rdr.
\eeq

Using the fact $u(0,t)=0$, it holds 
\beq\label{68}
(u(r))^2=2\displaystyle\int_0^r\frac{u(r)}{r}\partial_r(u(r))rdr
\le 2\left(\displaystyle\int_{\mathbb{R}^+}\left(\frac{u(r)}{r}\right)^2rdr\right)^\frac{1}{2}\left(\displaystyle\int_{\mathbb{R}^+}(\partial_ru)^2rdr\right)^\frac{1}{2}.
\eeq
Therefore, we obtain
$$\displaystyle\int_{\mathbb{R}^+}\frac{u^4}{r}dr
\le 2\displaystyle\int_{\mathbb{R}^+}\frac{u^2}{r}\,dr\cdot \left(\displaystyle\int_{\mathbb{R}^+}\left(\frac{u(r)}{r}\right)^2rdr\right)^\frac{1}{2}\left(\displaystyle\int_{\mathbb{R}^+}(\partial_ru)^2rdr\right)^\frac{1}{2}
\le 2(E_0(u,v))^2,$$ 
which gives
\[
|\mathcal{R}_\lambda(u)|
\le C\left(\displaystyle\int_{\mathbb{R}^+}\left(\partial_ru+\frac{u}{r}\right)^2\,rdr\right)^\frac{1}{2}\left(\displaystyle\int_{\mathbb{R}^+}\frac{u^4}{r^2}rdr\right)^\frac{1}{2}\\
\leq C \big(E_0(u,v)\big)^\frac{3}{2}.
\]
Putting this estimate into \eqref{est-E0R}, we obtain
$$E_0[u,v]\leq  C\epsilon^2+C\big(E_0[u,v]\big)^\frac{3}{2}.$$ 
By the property of the function  $f(x)=x-Cx^\frac{3}{2}=x(C-x^\frac{1}{2})$, we know that $E_0[u,v]$ is either around $0$ or around $C$. Since $E_0[u,v]$ is  continuous in $t$ and $E_0[u_0,v_0]\leq  \epsilon^2$ for sufficiently small $\epsilon$, it is not hard to see that 
\[\left.\begin{array}{l}
 C-E_0[u,v]^\frac{1}{2}>\frac{1}{2}
\end{array}\right.\]
for any $t\in[0,\tilde T_1]$ with $\tilde 0<T_1<T_1$.
Therefore $$C\epsilon^2\geq E_0[u]-E_0[u]^\frac{3}{2}=E_0[u](1-E_0[u]^\frac{1}{2})\ge \frac{1}{2}E_0[u]
$$
which implies
$$ E_0[u]\leq C \epsilon^2 \hbox{on }[0,\tilde T_1].$$
By standard continuity argument, we can show 
$$ E_0[u]\leq C \epsilon^2 \hbox{on }[0, T_1].$$

\paragraph{\bf Step 3.} Combining last two steps, we have shown the uniform estimate \eqref{u-est-alpha}, by which 
we can extend $\lambda(t)$ to $(0,T]$.

%

In order to show the estimate \eqref{cond3}, by the equation \eqref{sec4pf2} of $\lambda$, we obtain for sufficient small $\epsilon>0$
\beq\notag
\begin{split}
\left|\frac{\dot{\lambda}}{\lambda^2}\right|
=\left|\frac{\langle\phi_t,J_\lambda\rangle}{\alpha(\lambda,t)}\right|\lambda
=\left|\frac{\langle \phi_t,J_\lambda\rangle}{\alpha(\lambda(t),t)+\langle J,J\rangle-\langle J,J\rangle}\right|\lambda
\le\frac{|\langle\phi_t,J_\lambda\rangle|\lambda}{|\langle J,J\rangle|-C\epsilon}
\leq C\epsilon,
\end{split}
\eeq
where we have used the fact
\beq\notag
|\langle\phi_t,J_\lambda\rangle|\lambda
\le C\left(\int_{\mathbb{R}^+}(\phi_t)^2\,rdr\right)^\frac{1}{2}\left(\int_{\mathbb{R}^+}(J_\lambda)^2\,rdr\right)^\frac{1}{2}\lambda(t)
\leq C\epsilon \left(\int_{\mathbb{R}^+}(J)^2\,rdr\right)^\frac{1}{2}\leq C\epsilon.
\eeq
This completes the proof of proposition. 
\endpf


\medskip
Finally, we summarize some useful estimates for future use. By the decomposition \eqref{jte} and \eqref{phiiu},  it holds 
\[
\int_{\R^+}|\phi_t|^2r^{2+\delta}\,dr
 \leq  \frac12\int_{\R^+} \frac{\dot \lambda^2}{\lambda^2} J_\lambda^2 r^{2+\delta}\,dr+
 \frac12\int_{\R^+} u^2_t r^{2+\delta}\,dr,
\] 
where
\[
\int_{\R^+} \frac{\dot\lambda^2}{\lambda^2} J^2_\lambda r^{2+\delta}\,dr
=\frac{\dot\la^2}{\la^{5+\delta}} \int_{\R^+} J^2_\lambda (r\la)^{2+\delta}\,dr\la,
\]
which is bounded when $k\geq 4.$
Similarly, by \eqref{I_def2},
\[
\int_{\R^+}|\phi_r|^2r^{2+\delta}\,dr
 \leq  \frac12\int _{\R^+} J_\lambda^2 r^{\delta}\,dr+
 \frac12\int_{\R^+} u^2_r r^{2+\delta}\,dr,
\] 
where
\[
\int_{\R^+}  J^2_\lambda r^{\delta}\,dr
=\frac{1}{\la^{1+\delta}} \int_{\R^+} J^2_\lambda (r\la)^{\delta}\,dr\la.
\]

By \eqref{vint0}-\eqref{2930}, it is easy to check that initially
\[
\int_{\R^+}|u_t|^2(r,0) r^{2+\delta}\,dr\quad \hbox{or equivalently}\quad \int_{\R^+}|\phi_1(r)|^2 r^{2+\delta}\,dr,
\]
\[
\int_{\R^+}|u_r|^2(r,0) r^{2+\delta}\,dr\quad \hbox{or equivalently}\quad \int_{\R^+}|\phi_0'(r)|^2 r^{2+\delta}\,dr
\]
and
\[
\int_{\R^+}|h_t|^2(r,0) r^{2+\delta}\,dr
\]
are all bounded. Then \eqref{2.32} and \eqref{2.33} tells that
\beq\label{allests}
\int_{\R^+} (|u_t|^2+|u_r|^2+|\phi_t|^2+|\phi_r|^2+|h_t|^2)(r,t) r^{2+\delta}\,dr<C,
\qquad\int_0^T \int_{\R^+} |h_{tr}|^2 r^{2+\delta}\,drdt<C
\eeq
as long as the solution is smooth.


\section{More estimates on trajectory $\lambda(t)$} 
\label{sec5}
\setcounter{equation}{0}

Derivation of singularity formation heavily relies on the following second order equation of $\lambda$. 

\begin{proposition}\label{prop5.1}
Let $\lambda(t)$ and $u(r,t)$ be defined as in Proposition \ref{prop3.1}. Then as long as $\lambda(t)$ is smooth, it satisfies the following second order equation
\beq\label{key1}
C_0(\ddot{\lambda} -2\frac{\dot\la^2}{\la})=2\langle u_t,\dot{J}_\lambda\rangle \lambda^3+\langle u,\ddot{J}_\lambda\rangle\lambda^3+\langle \mathcal{N}(u),J_\lambda\rangle\lambda^3-C_0\dot\lambda-\langle u_t,J_\lambda\rangle\lambda^3-\langle h_t,J_\lambda\rangle\lambda^3
\eeq
with initial data
\beq\label{initial_2}
\lambda(0)=\lambda_0=\v^{-4},\quad \dot\lambda(0)=\frac{C_0}{\pi} \v^{-7}.
\eeq
\end{proposition}

\pf First notice that by \eqref{31} and \eqref{cond2}, it holds
\[
\langle \dot{I}_\lambda,J_\lambda\rangle=\langle\phi_t,J_\lambda\rangle-\langle u_t,J_\lambda \rangle=\langle\phi_t,J_\lambda\rangle+\langle u,\dot{J}_\lambda \rangle.
\]
Taking another time derivative, we get
\[
\frac{d}{dt}\langle \dot{I}_\lambda,J_\lambda\rangle=\langle\phi_{tt},J_\lambda\rangle+\langle\phi_t,\dot{J}_\lambda\rangle+\langle u_t,\dot{J}_\lambda\rangle+\langle u,\ddot{J}_\lambda\rangle.
\]
By \eqref{3.39},
$ \phi_{tt}=-u_t-H_\lambda u-\dot{I}_\lambda+\mathcal{N}(u)-h_t$, 
then we know that,
\[
\langle \ddot{I}_\lambda,J_\lambda\rangle=-\langle u_t,J_\lambda\rangle-\langle\dot{I}_\lambda,J_\lambda\rangle+\langle \mathcal{N}(u),J_\lambda\rangle-\langle h_t,J_\lambda\rangle+2\langle u_t,\dot{J}_\lambda\rangle+\langle u,\ddot{J}_\lambda\rangle,
 \]
 where we use $ \langle  H_\lambda u, J_\lambda\rangle=0$ by \eqref{hleq} and \eqref{aljl}.
Recalling $\dot{I}_\lambda=\frac{\dot{\lambda}}{\lambda}J_\lambda$ in \eqref{jte} and \eqref{C0}, we know
$$\langle\ddot{I}_\lambda,J_\lambda\rangle=(\frac{\ddot\lambda}{\lambda}-2\frac{{\dot \lambda}^2}{\lambda^2})\langle J_\lambda,J_\lambda\rangle=C_0\frac{1}{\lambda}\frac{d}{dt}(\frac{\dot{\lambda}}{\lambda^2}).$$
Again using $\langle u_t,J_\lambda\rangle=-\langle u,\dot{J_\lambda}\rangle$ by \eqref{cond2}, we have 
 \[
C_0\frac{1}{\lambda}\frac{d}{dt}(\frac{\dot{\lambda}}{\lambda^2})=2\partial_t\langle u,\dot{J}_\lambda\rangle-\langle u,\ddot{J}_\lambda\rangle+\langle \mathcal{N}(u),J_\lambda\rangle-\langle \dot{I}_\lambda,J_\lambda\rangle-\langle u_t,J_\lambda\rangle-\langle h_t,J_\lambda\rangle.\]
Hence
\beq\notag
C_0\frac{d}{dt}(\frac{\dot{\lambda}}{\lambda^2})=2\partial_t[\langle u,\dot{J}_\lambda\rangle\lambda]-2\langle u,\dot{J}_\lambda\rangle \dot{\lambda}-\langle u,\ddot{J}_\lambda\rangle\lambda+\langle \mathcal{N}(u),J_\lambda\rangle\lambda-\langle \dot{I}_\lambda,J_\lambda\rangle\lambda-\langle u_t,J_\lambda\rangle\lambda-\langle h_t,J_\lambda\rangle\lambda,
\eeq
or equivalently,
\beq\label{key11}
C_0\frac{d}{dt}(\frac{\dot{\lambda}}{\lambda^2})=2\partial_t[\langle u,\dot{J}_\lambda\rangle\lambda]-2\langle u,\dot{J}_\lambda\rangle \dot{\lambda}-\langle u,\ddot{J}_\lambda\rangle\lambda+\langle \mathcal{N}(u),J_\lambda\rangle\lambda-C_0\frac{\dot\lambda}{\lambda^2}-\langle u_t,J_\lambda\rangle\lambda-\langle h_t,J_\lambda\rangle\lambda.
\eeq
Therefore, we obtain the equation \eqref{key1}.

For the initial data, by the definition of \eqref{lambda1} and \eqref{C0}, it is easy to compute 
\beq\label{initial_1}
\lambda(0)=\lambda_0=\v^{-4}, \quad \dot\lambda(0)=\frac{C_0}{\pi} \v^{-7},
\eeq
which completes the proof of proposition.
\endpf

\bigskip

The main purpose of this section is to show the following estimate.

\begin{proposition}\label{prop5.2}
Let $\lambda(t)$ be solutions to the initial value problem \eqref{key1} and \eqref{initial_2}. When $t\in[0,\v^{\frac{1}{2}}]$ or before blowup of $\lambda$, it holds
\beq\label{key3}
\left|\ddot{\lambda}(t)-2\frac{\dot{\la}^2}{\lambda}(t)\right|\leq 
C_1 c_0^{\frac{1}{2}}\frac{\dot{\la}^2}{\lambda}(t)
+C_1 \la^2(t) \left(\sup_{0\leq s\leq t}\frac{\dot{\la}^4}{\la^7}(s)+c_0\v^2\right),
\eeq
for some positive constant $C_1$.
\end{proposition}

We will utilize the bootstrapping arguments to prove Proposition \ref{prop5.2}. To this end, we need to show several lemmas.

\begin{lemma}\label{lemma5.3}
If we assume that 
\beq\label{tech2}
|\ddot{\lambda}(t)-2\frac{\dot{\la}^2}{\lambda}(t)|\leq 
2C_1 c_0^{\frac{1}{2}}\frac{\dot{\la}^2}{\lambda}(t)
+2C_1  \la^2(t)\left(\sup_{0\leq s\leq t}\frac{\dot{\la}^4}{\la^7}(s)+c_0\v^2\right),
\eeq
for small enough $\v$, $c_0$ satisfying $\v^{\frac12}\leq c_0$. 
When $t\in[0,\v^{\frac{1}{2}}]$ or before blowup of $\lambda$, it holds
\beq\label{tech4}
\frac{d}{dt}(\frac{\dot{\la}^4}{\la^7})>0,
\eeq
and 
\beq\label{tech5}
\frac{d}{dt}\la>0.
\eeq
\end{lemma}

\begin{remark}\label{remark7.2}
As a direct result of this lemma and Proposition \ref{prop5.2}, we can prove \eqref{tech4}-\eqref{tech5} for any smooth solution before blowup or $t=1$ by a continuity argument, which implies
\beq\label{key0}
 \lambda(t)\geq O(\v^{-4}),\quad \frac{\dot{\la}^4}{\la^7}(t)\geq 1,
\eeq
for any time before blowup or $t=1$.
\end{remark}

\pf
Direct computation implies 
\beq\label{dl4l70}
\begin{split}
\frac{d}{dt}\left(\frac{\dot{\la}^4}{\lambda^7}(t)\right)
= 4\ddot{\la}\frac{\dot{\la}^3}{\la^7}-7\frac{\dot{\la}^5}{\la^8}
=4\frac{\dot{\la}^3}{\la^7}\left(\ddot{\la}-2\frac{\dot{\la}^2}{\la}\right)+\frac{\dot{\la}^5}{\la^8}.
\end{split}
\eeq
Hence, by assumption \eqref{tech2},
\beq\label{dl4l7}
\begin{split}
\frac{d}{dt}\left(\frac{\dot{\la}^4}{\lambda^7}(t)\right)
= \frac{\dot{\la}^5}{\lambda^8}(t)
+\hbox{O}\left(c_0^{\frac{1}{2}}\frac{\dot{\la}^2}{\lambda^3}(t)+\sup_{0\leq s\leq t}\frac{\dot{\la}^4}{\la^7}(s)+c_0\v^2 \right)\frac{\dot{\la}^3}{\lambda^5}(t).
\end{split}
\eeq
By \eqref{tech2}, it holds initially
\beq\label{sec7pf3}
\begin{split}
\frac{d}{dt}\left(\frac{\dot{\la}^4}{\lambda^7}(t)\right)\big|_{t=0}
\geq \frac{\dot{\la}^5}{\lambda^8}(0)
-C_1\left(c_0^{\frac{1}{2}}\frac{\dot{\la}^2}{\lambda^3}(0)+\frac{\dot{\la}^4}{\la^7}(0)+c_0\v^2 \right)\frac{\dot{\la}^3}{\lambda^5}(0),
\end{split}
\eeq
where we have used the fact $\lambda_0=\v^{-4}>0$ and $\dot{\la}(0)=C\v^{-7}>0$.
Hence by \eqref{dl4l7}, it is easy to see that 
\[
\frac{d}{dt}\left(\frac{\dot{\la}^4}{\lambda^7}(t)\right)\big|_{t=0}
\geq C\v^{-3}(1-c_0^{\frac{1}{2}}-\v^{2}-c_0\v^4)
\geq C^* \v^{-3}>0,
\]
for some positive constant $C^*$.
Now we prove that \eqref{tech4} holds
for any time by contradiction. If \eqref{tech4} fails to hold, then there exists the first time $\bar T>0$, such that  $\frac{d}{dt}(\frac{\dot{\la}^4}{\la^7})(\bar T)=0$. 
For  $t\in[0,\bar T]$, by the fact $\frac{d}{dt}(\frac{\dot{\la}^4}{\la^7})(t)\geq 0$, we know that
$$\frac{\dot{\la}^4}{\lambda^7}(t)\geq \frac{\dot{\la}^4}{\lambda^7}(0)=C>0.$$ 
Combining this estimate with the fact $\lambda(t)>0$, it holds 
$$\dot\lambda(t)> 0,$$
which implies $$\la(t)>\la_0=\v^{-4}.$$ 
Then it is not hard to see
$$
0<\left(\sup_{0\leq s\leq t}\frac{\dot{\la}^4}{\la^7}(s) \right)\frac{\dot{\la}^3}{\lambda^5}(t)\leq  \frac{\dot{\la}^4}{\la^7}(t)\frac{\dot{\la}^3}{\lambda^5}(t)
\leq C\v^2\frac{\dot{\la}^5}{\lambda^8}(t),
$$
$$
0<\frac{\dot{\la}^3}{\lambda^5}(t)\leq  C\frac{\dot{\la}^3}{\lambda^5}(t) \left(\frac{\dot{\la}^4}{\lambda^7}(t)\right)^{\frac{1}{2}}\left(\v^4\la\right)^\frac{1}{2}\leq C\v^2 \frac{\dot{\la}^5}{\lambda^8}(t) .
$$
Combining all these facts with \eqref{dl4l7}, we prove that  $\frac{d}{dt}(\frac{\dot{\la}^4}{\la^7})(T^*)>0$. Hence we find a contradiction, which proves the estimate \eqref{tech4}. 
By \eqref{tech4} and the above arguments we can also prove \eqref{tech5}, which completes the proof of lemma. 
\endpf

\begin{lemma}\label{lemma5.5} Let $\lambda(t)$ satisfy the same bootstrapping assumption in Lemma \ref{lemma5.3}. 
When $t\in[0,\v^{\frac{1}{2}}]$ or before blowup of $\lambda$, it holds
\beq\label{tech3}
\left|\ddot{\lambda}(t)-2\frac{\dot{\la}^2}{\lambda}(t)\right|\leq 
C_1 c_0^{\frac{1}{2}}\frac{\dot{\la}^2}{\lambda}(t)
+C_1 \la^2(t) \left(\sup_{0\leq s\leq t}\frac{\dot{\la}^4}{\la^7}(s)+c_0\v^2\right).
\eeq
\end{lemma}
This is the main technical estimate of our paper, which will be proved in Section \ref{sec6} below.

\bigskip
\noindent{\bf Proof of Proposition \ref{prop5.2}}. \quad
By the bootstrapping argument, Lemma \ref{lemma5.3} and Lemma \ref{lemma5.5}, we only need to show the following estimate initially
\beq\label{sec5pf2}
\left|\ddot{\lambda}(0)-2\frac{\dot{\la}^2}{\lambda}(0)\right|
\leq 
C_1 c_0^{\frac{1}{2}}\frac{\dot{\la}^2}{\lambda}(0)
+C_1 \la^2(0)\left(\frac{\dot{\la}^4}{\la^7}(0)+c_0\v^2\right).
\eeq
In fact, we need to estimate all the terms on the right hand side of \eqref{key1} at initial time. We first deal with all terms including $u$. It holds
\beq\label{p5.2pf1}
\begin{split}
\left|\langle u_0,\ddot{J}_{\lambda_0}\rangle\lambda_0^3\right|
\leq C\int_{\R^+}\left|u_0\right|\,\left(
\left|\ddot{\la}(0)\right|\la_0^2+\dot{\la}^2(0)\la_0\right)F^{4,4}_{\la_0}\,rdr.
\end{split}
\eeq
For the first term, we can estimate similar as \eqref{sec4pf3}
\beq\notag
\begin{split}
 \left|\ddot{\la}(0)\right|\int_{\R^+}\left|u_0\right|
\la_0^2F^{4,4}_{\la_0}\,rdr
= \left|\ddot{\la}(0)\right|\int_{\R^+}\left|u_0\left(\frac{r}{\la_0}\right)\right|\,F^{4,4}\,rdr
\leq c_0\v\left|\ddot{\la}(0)\right|.
\end{split}
\eeq
For the second term, it holds
\beq\notag
\begin{split}
\int_{\R^+}\left|u_0\right|\dot{\la}^2(0)
\la_0F^{4,4}_{\la_0}\,rdr
\leq \frac{\dot\la^2}{\la}(0)\left(\int_{\R^+}\frac{\left|u_0\right|^2}{r^2}\,rdr\right)^{\frac12}\left(\int_{\R^+}F^{8,8}_{\la_0}\la_0^4\,r^3dr\right)^{\frac12}
\leq c_0\v\frac{\dot\la^2}{\la}(0).
\end{split}
\eeq
Putting these into \eqref{p5.2pf1}, we have
\beq\label{p5.2pf2}
\begin{split}
\left|\langle u_0,\ddot{J}_{\lambda_0}\rangle\lambda_0^3\right|
\leq Cc_0\v\left(\left|\ddot{\la}(0)-2\frac{\dot\la^2}{\la}(0)\right|+\frac{\dot\la^2}{\la}(0)\right).
\end{split}
\eeq
By the definition of $\mathcal{N}(u)$, it is not hard to conclude 
\beq\label{p5.2pf3}
\left|\langle \mathcal{N}(u_0),J_{\lambda_0}\rangle\lambda_0^3\right|
\leq C\la_0^3\int_{\R^+}\frac{|u_0|^2}{r^2}F^{4,4}_{\la_0}\,rdr
\leq Cc_0^2\v^2\la_0^3\leq Cc_0^2\frac{\dot\la^2}{\la}(0).
\eeq

For all the terms including $u_t$,  by the definition of $u$, we have
\beq\notag
u_t(r,0)=\phi_1(r)-\frac{\dot{\la}}{\la}(0)J_{\la_0}=g_0(r).
\eeq
Similar computations as \eqref{p5.2pf2}-\eqref{p5.2pf3} imply
\beq\label{p5.2pf4}
\left|\langle u_t(r,0),\dot{J}_{\lambda_0}\rangle \lambda_0^3\right|+\left|\langle u_t(r,0),J_{\lambda_0}\rangle\lambda_0^3\right|
\leq Cc_0\v\frac{\dot\la^2}{\la}(0).
\eeq
For other terms, it is not hard to see
\beq\label{p5.2pf5}
|\dot\la(0)|=\frac{C_0}{\pi}\v^{-7}\leq Cc_0^2\v^2\frac{\dot\la^2}{\la}(0).
\eeq
By the definition \eqref{eqnv}, we know $h_t(0)=\phi_t(0)=\phi_1=g_0+\la_1J_{\la_0}$. Hence
\beq\label{p5.2pf6}
\left|\langle h_t(r,0),J_{\lambda_0}\rangle \lambda_0^3\right|
\leq Cc_0^2\v^2\frac{\dot\la^2}{\la}(0)+C_0\la_0\la_1
\leq Cc_0^2\v^2\frac{\dot\la^2}{\la}(0).
\eeq
Combining \eqref{p5.2pf2}-\eqref{p5.2pf6} with equation \eqref{key1} for sufficiently small $\v$, we conclude the estimate \eqref{sec5pf2} at initial time.
Therefore, the Proposition \ref{prop5.2} can be proved by using bootstrapping argument and Lemma \ref{lemma5.3}-Lemma \ref{lemma5.5}.
\endpf

%
%
%
%
%


\section{Proof of bootstrapping Lemma \ref{lemma5.5}}
\label{sec6}
\setcounter{equation}{0}

We need the following Morawetz estimate which has been proved as Proposition 6.1 in \cite{RS}.

\begin{proposition}
\label{prop6.1}
Let $\psi$ be a smooth function on $[t_0,t_1]\times(0,\infty)$ with the following uniform bounds for all $0\leq r\leq 1$
\beq\label{RS119}
|\psi|\leq C(\psi(t))r,\quad |\psi_t|+|\psi_r|\leq C(\psi)
\eeq
while decaying sufficiently rapidly at $r=\infty$. Suppose that 
\beq\label{RS120}
\partial^2_t\psi+\tilde H_\la \psi=\partial_t \mathcal G+\mathcal H,
\eeq
while $\dot{\lambda}\geq0$, and $\dot{\lambda}\lambda^{-2}\leq C\v$ for all times $t_0\leq s\leq t_1$, then the following estimate is valid with the constant independent of $\lambda$ and $\delta$
\beq\label{RS121} 
\begin{split}
E_\delta[\psi](t_0,t_1)
\leq & C\delta^{-1}
\left[
\int_{t_0}^{t_1}\int_{\R^+}\lambda^{-1}(\lambda r)^{\delta}\left[(\partial_r\mathcal G)^2+\v^2(\lambda \mathcal G)^2+\mathcal H^2\right](s)\,r^2drds\right.\\
&\left.+\sup\limits_{t_0\leq s\leq t_1}\int_{\R^+}\lambda^{-1}\frac{(\lambda r)^{\delta}}{1+r^{\delta}}\mathcal G^2(s)\,rdr+E_{\delta}[\psi](t_0,t_0)\right],
\end{split}
\eeq
where 
\beq\label{RS118}
\begin{split}
E_\delta[\psi](t_0,t_1)=&\sup_{t_0\leq s\leq t_1}\int_{\mathbb{R^+}}\la^{-1}\frac{(\la r)^\delta}{1+r^\delta}
\left[(L\psi)^2+\frac{\psi^2}{r^2}\right](s)\,r\,dr\\
&+ \int_{t_0}^{t_1}\int_{\mathbb{R^+}}\left[\la^{-1}\frac{(\la r)^\delta}{(1+r^\delta)^2r}
(L\psi)^2+\frac{(\la r)^\delta}{1+r^\delta}\frac{\psi^2}{r^3}\right](s)\,r\,drds,
\end{split}
\eeq
and $L=\partial_r+\partial_t$ is the outgoing null derivatives.
\end{proposition}

The goal is to first prove that 
\begin{lemma}\label{prop_tech1}
Let $u=w_0+w$ be the decomposition defined in \eqref{uw0w}. Let us assume that 
\beq\label{tech1}
|\ddot{\lambda}(t)-2\frac{\dot{\la}^2}{\lambda}(t)|\leq 
2C_1 c_0^{\frac{1}{2}}\frac{\dot{\la}^2}{\lambda}(t)
+2C_1  \la^2(t)\left(\sup_{0\leq s\leq t}\frac{\dot{\la}^4}{\la^7}(s)+c_0\v^2 \right),
\eeq
for some positive constant $C_1$ and small constant $c_0\geq \v^{\frac{1}{2}}$, 
and $\lambda(t)$ is monotonically nondecreasing. For any time interval $[t_0,t_1]$ in $[0,1]$ and before the blowup of $\la$, we assume 
\beq\label{extassp}
\int_{t_0}^{t_1}\frac{\dot{\lambda}^4}{\lambda^7}(s)\,ds\leq \v
\eeq
for some small enough $\v>0$. Then one has the following estimate
\beq\label{tech21}
E_\delta[A_\lambda w](t_0,t_1)
\leq  C \delta^{-1}
\left(
E_\delta(A_\lambda w)(t_0,t_0)+
\v^3+\v\sup\limits_{t_0\leq t\leq t_1}\frac{\dot{\lambda}^4}{\lambda^7}(s)
\right).
\eeq
Here $\delta\ll 1$ is a positive constant. The constant $C$ depends on $C_1$ in \eqref{tech1}  but is independent of 
$\v$, $c_0$ and $\delta$.

\end{lemma}

%

\pf By the equation  \eqref{3.39} and the fact $u=w_0+w$, we have
\beq\label{eqnw}
\partial_{tt} w+w_t+H_\lambda w=-\ddot{I}_\lambda-\dot{I}_\lambda+\mathcal{N}(u)-h_t-\partial_{tt}w_0-\partial_tw_0-H_\lambda w_0.
\eeq
Operating $A_\lambda$ on \eqref{eqnw}, it is easy to get 
\beq\label{134}
\begin{split}
&\partial_{tt} (A_\lambda w)+(A_\lambda w)_t+\tilde H_\la A_\lambda w\\
=&\partial_tA_\la\partial_tw_0
-\partial_t(A_\la\partial_tw_0)+2\partial_t(\partial_tA_\la w)
-\partial_{tt}A_\la w+A_\lambda \mathcal{N}(u)-A_\lambda  h_t\\
& -(A_\lambda w_0)_t+\partial_tA_\lambda w-\partial_tA_\lambda w_0
\end{split}
\eeq
where we use the fact 
$$
A_\lambda H_\lambda=\tilde H_\lambda A_\lambda
$$
and the equations \eqref{w0AH} and \eqref{AId1}.
To deal with the damping term, we denote 
\beq\label{defWa}
W=e^{\frac{t}{2}}A_\lambda w.
\eeq
By direct computation and using \eqref{134}, it holds
\beq\label{134damp}
\begin{split}
&\partial_{tt} (W)+\tilde H_\la W=\partial_t \mathcal G+\mathcal H \\
\end{split}
\eeq
where
\beq\label{Gdef}
\begin{split}
\mathcal G=-e^{\frac{t}{2}}A_\la\partial_t w_0
+2e^{\frac{t}{2}}\partial_tA_\la w:=\Phi_1+\Phi_2
\end{split}
\eeq
\beq\label{Hdef}
\begin{split}
\mathcal H=e^{\frac{t}{2}}
\left[\partial_tA_\la\partial_tw_0-\partial_{tt}A_\la w+A_\la\mathcal N(u)-\frac{1}{2}A_\la(\partial_tw_0)-A_\la h_t+\frac14A_\la w\right]:=\sum\limits_{i=1}^6 \Psi_i.
\end{split}
\eeq

By the definition of $w$ in \eqref{uw0w}, it is not hard to verify that $W=e^{\frac{t}{2}}A_\lambda w$ satisfies the assumption \eqref{RS119} for any $0\leq r\leq 1$. Hence, we may use Proposition \ref{prop6.1} to control $E_\delta[W](t_0,t_1)$ by a bound as \eqref{RS121}, with $\mathcal G$ and $\mathcal H$ in \eqref{Gdef} and \eqref{Hdef}, respectively. More precisely
\beq\label{sec6pf3} 
\begin{split}
E_\delta[W](t_0,t_1)
\leq & C\delta^{-1}
\left[
\int_{t_0}^{t_1}\int_{\R^+}\lambda^{-1}(\lambda r)^{\delta}\left[(\partial_r\mathcal G)^2+\v^2(\lambda \mathcal G)^2+\mathcal H^2\right](s)\,r^2drds\right.\\
&\left.+\sup\limits_{t_0\leq s\leq t_1}\int_{\R^+}\lambda^{-1}\frac{(\lambda r)^{\delta}}{1+r^{\delta}}\mathcal G^2(s)\,rdr+E_{\delta}[W](t_0,t_0)\right],
\end{split}
\eeq
where $E_\delta[W](t_0,t_1)$ is given as in \eqref{RS118}.
%

 The estimates corresponding to $\mathcal G$ and the first three terms of $\mathcal H$ are exactly the same as in \cite{RS} (cf. estimates (135)-(143) in \cite{RS}). For readers' convenience, we list these estimates below.
\beq\label{136-138}
\begin{split}
\int_{t_0}^{t_1}\int_{\mathbb{R^+}}
\frac{(\lambda r)^\delta}{\lambda}\left(|\partial_r \Phi_1|^2+|\lambda\Phi_1|^2\right)r^2drds
+\sup_{t_0\leq t\leq t_1}\int_{\mathbb{R^+}}\frac{1}{\lambda}\frac{(\lambda r)^\delta}{1+r^\delta}| \Phi_1|^2rdr
\leq C\v^3+C\epsilon\sup_{t_0\leq t\leq t_1}\frac{\dot\lambda^4}{\lambda^7}(t),
\end{split}
\eeq
%
\beq\label{139-141}
\begin{split}
\int_{t_0}^{t_1}\int_{\mathbb{R^+}}\frac{(\lambda r)^\delta}{\lambda}
\left(|\partial_r \Phi_2|^2+|\lambda \Phi_2|^2\right)r^2drds
+\sup_{t_0\leq t\leq t_1}\int_{\mathbb{R^+}}\frac{1}{\lambda}\frac{(\lambda r)^\delta}{1+r^\delta}| \Phi_2|^2rdr
\leq C\epsilon^2E_\delta[A_\lambda w](t_0,t_1),
\end{split}
\eeq
%
\beq\label{135}
\begin{split}
\int_{t_0}^{t_1}\int_{\mathbb{R^+}}\frac{(\lambda r)^\delta}{\lambda}| \Psi_1|^2r^2drds
\leq C\epsilon\sup_{t_0\leq t\leq t_1}\frac{\dot\lambda^4}{\lambda^7}(t),
\end{split}
\eeq
\beq\label{142}
\begin{split}
\int_{t_0}^{t_1}\int_{\mathbb{R^+}}\frac{(\lambda r)^\delta}{\lambda}| \Psi_2|^2r^2drds
\leq C\epsilon^2E_\delta[A_\lambda w](t_0,t_1),
\end{split}
\eeq
\beq\label{143}
\begin{split}
\int_{t_0}^{t_1}\int_{\mathbb{R^+}}\frac{(\lambda r)^\delta}{\lambda}| \Psi_3|^2r^2drds
\leq C\epsilon\sup_{t_0\leq t\leq t_1}\frac{\dot\lambda^4}{\lambda^7}(t)+t_1^\delta\epsilon^2E_\delta[A_\lambda w](t_0,t_1).
\end{split}
\eeq

To estimate $\Psi_4$, notice that it is similar to the term $\Phi_1$. Since $\lambda_0=\v^{-4}$ and $\lambda$ is nondecreasing, we have 
\beq\label{bdlambda2}
\lambda(t)\geq \lambda_0=\v^{-4}.
\eeq
Combining this fact with the estimate \eqref{136-138}, we obtain
\beq\label{estpsi4}
\int_{t_0}^{t_1}\int_{\mathbb{R^+}}
\frac{(\lambda r)^\delta}{\lambda}|\Psi_4|^2r^2drds
\leq
\v^8\int_{t_0}^{t_1}\int_{\mathbb{R^+}}
\frac{(\lambda r)^\delta}{\lambda}|\lambda\Phi_1|^2r^2drds
\leq C\epsilon\sup_{t_0\leq t\leq t_1}\frac{\dot\lambda^4}{\lambda^7}(t).
\eeq


To estimate $\Psi_5$, by the definition of $A_\lambda$, it holds
\beq\label{sec6pf1}
\begin{split}
\int_{t_0}^{t_1}\int_{\mathbb{R}^+}\frac{1}{\lambda}(\lambda r)^\delta|\Psi_5|^2r^2drds
\leq C\int_{t_0}^{t_1}\lambda^{\delta-1}\int_{\mathbb{R}^+}\left(|h_{tr}|^2+\frac{|h_t|^2}{r^2}\right)r^{2+\delta}drds.
\end{split}
\eeq
By the energy estimate \eqref{2.33} and \eqref{bdlambda2}, we have
\beq\notag
\int_{t_0}^{t_1}\lambda^{\delta-1}\int_{\mathbb{R}^+}|h_{tr}|^2r^{2+\delta}drds
\leq C\sup(\lambda^{\delta-1})\leq C\v^{4-4\delta}\leq C\v^3.
\eeq
By the energy estimates \eqref{2.17} and \eqref{2.33}, we have the following estimate for any $t\in[0,1]$
\beq\notag
\int_{\mathbb{R}^+}\frac{|h_t|^2}{r^2}r^{2+\delta}dr
=\int_0^1\frac{|h_t|^2}{r^2}r^{2+\delta}\,dr+\int_1^\infty\frac{|h_t|^2}{r^2}r^{2+\delta}\,dr
\leq C.
\eeq
Combining this fact with \eqref{bdlambda2}, it holds
\beq\notag
\int_{t_0}^{t_1}\lambda^{\delta-1}\int_{\mathbb{R}^+}\frac{|h_t|^2}{r^2}r^{2+\delta}drds
\leq C\sup(\lambda^{\delta-1})\leq C\v^{4-4\delta}\leq C\v^3.
\eeq
Putting all these estimates into \eqref{sec6pf1}, we obtain the required estimate for $\Psi_5$
\beq\label{sec6pf2}
\begin{split}
\int_{t_0}^{t_1}\int_{\mathbb{R}^+}\frac{1}{\lambda}(\lambda r)^\delta|\Psi_5|^2r^2drds
\leq C\v^3.
\end{split}
\eeq
To estimate $\Psi_6$, by the definition of $A_\lambda$ and \eqref{uw0w}, we have
\beq\label{sec6pf4}
\begin{split}
\int_{t_0}^{t_1}\int_{\mathbb{R}^+}\frac{1}{\lambda}(\lambda r)^\delta|\Psi_6|^2 r^2drds
\leq C\int_{t_0}^{t_1}\lambda^{\delta-1}\int_{\mathbb{R}^+}\left[|A_\lambda u|^2+|A_\lambda w_0|^2\right]r^{2+\delta}drds.
\end{split}
\eeq
For the first term, by the estimate \eqref{allests} and decay estimate \eqref{147} in Appendix Lemma \ref{decayu}, we have
\beq\label{sec6pf5}
\begin{split}
&\int_{t_0}^{t_1}\lambda^{\delta-1}\int_{\mathbb{R}^+}|A_\lambda u|^2r^{2+\delta}drds\\
\leq &C\int_{t_0}^{t_1}\lambda^{\delta-1}\int_{\mathbb{R}^+}\left[|\partial_r u|^2+\frac{|u|^2}{r^2}\right]r^{2+\delta}drds\\
\leq &C\int_{t_0}^{t_1}\lambda^{\delta-1}\left(C+\int_{\{r\leq 3s\}}|u|^2r^{\delta}dr+\int_{\{r\geq 3s\}}|u|^2r^{\delta}dr\right)ds \\
\leq & C\v^3.
\end{split}
\eeq
For the second term involving $w_0$, by the estimate \eqref{58a3}, we have 
\beq\label{sec6pf6}
\begin{split}
&\int_{t_0}^{t_1}\lambda^{\delta-1}\int_{\mathbb{R}^+}|A_\lambda w_0|^2r^{2+\delta}drds\\
\leq &C\int_{t_0}^{t_1}\dot{\lambda}^4\lambda^{\delta-7}\int_{\mathbb{R}^+}\left[F^{3,5}_\lambda+F^{5,3}_\lambda\right]^2r^{2+\delta}drds\\
\leq &C \int_{t_0}^{t_1}\dot{\lambda}^4\lambda^{-10}\,ds\\
\leq & C\sup(\lambda^{-3}) \int_{t_0}^{t_1}\dot{\lambda}^4\lambda^{-7}\,ds\leq C\v^3,
\end{split}
\eeq
where we have used the assumption \eqref{extassp} in last step.
Putting \eqref{sec6pf5} and \eqref{sec6pf6} into \eqref{sec6pf4}, we obtain
\beq\label{sec6pf7}
\begin{split}
\int_{t_0}^{t_1}\int_{\mathbb{R}^+}\frac{1}{\lambda}(\lambda r)^\delta|\Psi_6|^2 r^2drds\leq C\v^3.
\end{split}
\eeq
Therefore combining all estimates \eqref{136-138}-\eqref{143}, \eqref{estpsi4}, \eqref{sec6pf2}, \eqref{sec6pf7} with \eqref{sec6pf3}, we conclude 
\beq\label{sec6pf8} 
\begin{split}
E_\delta[A_\lambda w](t_0,t_1)\leq CE_\delta[W](t_0,t_1)
\leq  C \delta^{-1}
\left(
E_\delta(W)(t_0,t_0)+
\v^3+\v\sup\limits_{t_0\leq t\leq t_1}\frac{\dot{\lambda}^4}{\lambda^7}(s)
\right),
\end{split}
\eeq
where we have used the definition \eqref{defWa} of $W$ in first inequality. 
Comparing $E_\delta(W)(t_0,t_0)$ with $E_\delta(A_\lambda w)(t_0,t_0)$, we only need to estimate the following term at $t_0$
\beq\label{sec6pf9}
\begin{split}
\int_{\mathbb{R^+}}\la^{-1}\frac{(\la r)^\delta}{1+r^\delta}
|A_\lambda w|^2\,r\,dr\leq \lambda^{\delta-1}\int_{\mathbb{R}^+}\left[|A_\lambda u|^2+|A_\lambda w_0|^2\right]\frac{r^{1+\delta}}{1+r^{\delta}}dr.
\end{split}
\eeq
For the first term, by the definition of $A_\lambda$, we have
\beq\label{sec6pf10}
\begin{split}
\lambda^{\delta-1}\int_{\mathbb{R}^+}|A_\lambda u|^2\frac{r^{1+\delta}}{1+r^{\delta}}dr
\leq C\lambda^{\delta-1}\int_{\mathbb{R}^+}\left[|\partial_r u|^2+\frac{|u|^2}{r^2}\right]\,rdr
\leq  C\v^3,
\end{split}
\eeq
where we have used \eqref{cond1} of Proposition \ref{prop3.1} in last inequality. 
For the second term involving $w_0$, by the estimate \eqref{58a3}, we have 
\beq\label{sec6pf11}
\begin{split}
\lambda^{\delta-1}\int_{\mathbb{R}^+}|A_\lambda w_0|^2\frac{r^{1+\delta}}{1+r^{\delta}}dr
\leq C\dot{\lambda}^4\lambda^{\delta-7}\int_{\mathbb{R}^+}\left[F^{3,5}_\lambda+F^{5,3}_\lambda\right]^2r^{1+\delta}dr
\leq  C\v\sup\limits_{t_0\leq t\leq t_1}\frac{\dot{\lambda}^4}{\lambda^7}(s).
\end{split}
\eeq
Combining \eqref{sec6pf8}-\eqref{sec6pf11}, we conclude \eqref{tech21}, which completes the proof of proposition.
\endpf

\bigskip
Now using Lemma \ref{prop_tech1} and bootstrapping argument with \eqref{tech1} (For extension of the interval from $t_0$ to $t_1$ to $0$ to $1$, we directly apply the Proposition $6.3$ in \cite{RS}), on any interval before blowup, we have the following result.

\begin{lemma}\label{prop_tech2}
Let $u=w_0+w$ be the decomposition defined in \eqref{uw0w}. Assume that \eqref{tech1} is valid and $\lambda(t)$ is monotonically nondecreasing. One has the following estimate for $t\in[0,1]$ or before the blowup of $\la$,
\beq\label{tech22}
E_\delta[A_\lambda w](0,t)
\leq C \delta^{-1}\left(c_0^2\v^2+\v\sup_{0\leq s\leq t}\frac{\dot{\la}^4}{\la^7}(s)\right)
\eeq
where $\delta\ll 1$ is positive constant. The constant $C$ depends on $C_1$ in \eqref{tech1},  but is independent of 
$\v$, $c_0$ and $\delta$. In particular, for any fixed $0<\delta\ll 1$, we have
\beq\label{tech23}
\int_{\mathbb{R}^+}
\lambda^{-1}\frac{(\lambda r)^{\delta}}{1+r^{\delta}}
\left[(LA_\lambda w)^2+\frac{(A_\lambda w)^2}{r^2}\right](t)
\,r\, dr
\leq C\left(c_0^2\v^2 \v\sup_{0\leq s\leq t}\frac{\dot{\la}^4}{\la^7}(s)\right).
\eeq
\beq\label{tech24}
\int_{\mathbb{R}^+}\lambda^{-1}\frac{(\la r)^\delta}{1+r^\delta}\frac{w^2}{r^4}\,r\, dr
\leq C\left(c_0^2\v^2+\v\sup_{0\leq s\leq t}\frac{\dot{\la}^4}{\la^7}(s)\right).
\eeq
\end{lemma}

The proof of this lemma is same as those in \cite{RS} (cf. proof of estimate (128) on Page 227), so we omit it.

\bigskip
\noindent{\bf Proof of Lemma \ref{lemma5.5}}. \quad 
Recall that the second order equation of $\la$ is given by \eqref{key1} as follows
\[
C_0(\ddot{\lambda} -2\frac{\dot\la^2}{\la})=2\langle u_t,\dot{J}_\lambda\rangle \lambda^3+\langle u,\ddot{J}_\lambda\rangle\lambda^3+\langle \mathcal{N}(u),J_\lambda\rangle\lambda^3-C_0\dot\lambda-\langle u_t,J_\lambda\rangle\lambda^3-\langle h_t,J_\lambda\rangle\lambda^3.
\]
We need to estimate each term of the right hand side under the bootstrapping assumption \eqref{tech2}.
The proof of first three terms exactly follow the proof in Lemma 5.5 in \cite{RS}
\beq\label{6001}
\left|\langle u_t,\dot{J}_\lambda\rangle \lambda^3\right|
\leq 
 Cc_0^{\frac{3}{4}}\frac{\dot{\la}^2}{\lambda}(t)
+Cc_0^\frac{1}{4}  \la^2(t)\left(c_0\v^2+\sup_{0\leq s\leq t}\frac{\dot{\la}^4}{\la^7}(s)\right),
\eeq
\beq\label{6002}
\left|\langle u,\ddot{J}_\lambda\rangle\lambda^3\right|
=\eta_1(t)\ddot \lambda+\eta_2(t)\frac{\dot \la^2}{\la}\quad\hbox{with}\quad
|\eta_1|\leq C \v,\quad |\eta_2|\leq C \v,
\eeq
\beq\label{6003}
\left|\langle \mathcal{N}(u),J_\lambda\rangle\lambda^3\right|
\leq C\v^2\frac{\dot{\la}^2}{\lambda}(t)
+\la^2(t)\left(c^2\v^2+\v\sup_{0\leq s\leq t}\frac{\dot{\la}^4}{\la^7}(s)\right).
\eeq
By \eqref{cond3} in Proposition \ref{prop3.1}, it is easy to see that
\beq\label{6.1.1}
|\dot\lambda|\leq C\v\la^2(t)\leq C\v \lambda^2(t)\sup_{0\leq s\leq t}\frac{\dot{\la}^4}{\la^7}(s),
\eeq
where we have used the fact that $\frac{\dot\lambda^4}{\lambda^7}$ is increasing in time by Lemma \ref{lemma5.3}. By the definition \eqref{phiiu} of $u$, we know $u_t=\phi_t-\dot I_\la=\phi_t-\frac{\dot\la}{\la}J_\la$. Thus, it holds
\beq\label{6.1.2}
|\langle u_t,J_\lambda\rangle\lambda^3|
\leq
\la^3\left(\int_{\R^+} |\phi_t|^2 rdr\right)^\frac{1}{2}\left(\int_{\R^+} |J_\la|^2 rdr\right)^\frac{1}{2}+C_0\dot\la
\leq
C\v\la^2(t) \sup_{0\leq s\leq t}\frac{\dot{\la}^4}{\la^7}(s),
\eeq
where we have used \eqref{cond3} in Proposition \ref{prop3.1} and \eqref{6.1.1}. For the last term with $h$, by H$\ddot{\hbox{o}}$lder's inequality, we have 
\beq\notag
\left|\langle h_t,J_\lambda\rangle\right|\la^3
\leq
C\la^3\left(\int_{\R^+} |h_t|^2 rdr\right)^\frac{1}{2}\left(\int_{\R^+} |J_\la|^2\,rdr\right)^\frac{1}{2}
=CC_0\la^2\left(\int_{\R^+} |h_t|^2 rdr\right)^\frac{1}{2}.
\eeq
By the energy estimate \eqref{2.17} and the fact $h_0\equiv0$, we have on time interval $[0,\v^{\frac12}]$
\beq\notag
\left|\langle h_t,J_\lambda\rangle\right|\la^3
\leq
C\la^2\left(\int_{\R^+} |h_t(r,0)|^2 rdr\right)^\frac{1}{2}+C\v^{\frac14}\la^2.
\eeq
By the equation \eqref{eqnv} of $h$, we know $h_t(r,0)=\phi_1=g_0+\la_1J_{\la_0}$. Thus
\beq\notag
\int_{\R^+} |h_t(r,0)|^2 \,rdr\leq Cc_0\v^2+C\frac{\la_1}{\la_0}\leq C\v.
\eeq
Combining these estimates together, we obtain
\beq\label{6.1.3}
\left|\langle h_t,J_\lambda\rangle \lambda^3\right|
\leq C \v^{\frac14}\la^2(t) \sup_{0\leq s\leq t}\frac{\dot{\la}^4}{\la^7}(s).
\eeq
Putting all estimates \eqref{6001}-\eqref{6.1.3} together and choosing small enough $c_0$ and $\v$, we can prove the bootstrapping Lemma \ref{lemma5.5}.
\endpf

\begin{remark}
In Section \ref{sec5}, we have shown that $\frac{\dot\lambda^4}{\lambda^7}$ is nondecreasing, so that by the computations in Section \ref{sec5}, we obtain that $\epsilon$ is much smaller than $\sup\frac{\dot\lambda^4}{\lambda^7}$. Therefore, we are going to delete $\epsilon$ terms in the future for simplicity.
\end{remark}


\section{Riccati equation and blowup}
\setcounter{equation}{0}

To obtain the singularity, we need to rewrite the equation \eqref{63} of $\lambda$ in terms of $u$ by the computation in Section \ref{sec5}.

\begin{lemma}\label{lemma6.1} Let $\lambda(t)$ be the solution to \eqref{63} obtained in Proposition \ref{prop3.1}. Then it holds
\beq\label{8.22}\begin{split}
&\left(C_0-2\langle u,r\partial_r J\rangle\right)\dot\lambda\\
=&\left(C_0\frac{\dot\lambda(0)}{\lambda_0^2}-2\dot\lambda_0\langle u_0,(r\partial_r J)_{\lambda_0}\rangle-\frac{2C_0}{\lambda_0}\right)\lambda^2+2C_0\lambda\\
&+\lambda^2\int_0^t\left\{-2\langle u,\dot{J}_\lambda\rangle\dot{\lambda}-\langle u,\ddot{J}_\lambda\rangle\lambda+\langle \mathcal{N}(u),J_\lambda\rangle\lambda \right\}ds\\
&+\lambda^2\int_0^t\left\{-\langle u_t,J_\lambda\rangle\lambda-\langle h_t,J_\lambda\rangle\lambda \right\}ds,
\end{split}
\eeq
where $C_0:=\langle J, J\rangle$. The coefficient of first term on the right hand side satisfies the following estimate
\beq\label{initial_3}
\frac{C_0\dot\lambda(0)}{\lambda_0^2}-2\dot\lambda(0)\langle u_0,(r\partial_r J)_{\lambda_0}\rangle-2\frac{C_0}{\lambda_0}\geq C\v
\eeq
for some constant $C> 0$ independent of $\v$. Furthermore, the coefficient on  the left hand side satisfies
\beq\label{sec5pf1}
\left|\langle u,r\partial_r J\rangle\right|\leq C\v.
\eeq
\end{lemma}

\pf 
Integrating the equation \eqref{key11}, we have
\beq\label{key2}
\begin{split}
&C_0\frac{\dot{\lambda}(t)}{\lambda(t)^2}-2\langle u(t),\dot{J}_\lambda\rangle\lambda \\
=&C_0\frac{\dot{\lambda}(0)}{\lambda_0^2}-2\dot\lambda(0)\langle u_0,(r\partial_rJ)_{\lambda_0}\rangle+2C_0\left(\frac{1}{\lambda(t)}-\frac{1}{\lambda_0}\right)\\
&+\int_0^t\left\{-2\langle u,\dot{J}_\lambda\rangle\dot{\lambda}-\langle u,\ddot{J}_\lambda\rangle\lambda+\langle \mathcal{N}(u),J_\lambda\rangle\lambda\right\} ds+\displaystyle\int_0^t\left\{-\langle u_t,J_\lambda\rangle\lambda-\langle h_t,J_\lambda\rangle\lambda \right\} ds
\end{split}
\eeq
where we use
$$
-2\langle u,\dot J_\lambda\rangle\lambda=-2\dot\lambda\langle u,(r\partial_r J)_\lambda\rangle
$$
because of 
$
\dot J_\lambda=\frac{\dot\lambda}{\lambda}(r\partial_r J)_\lambda
$
in \eqref{jte}. Therefore, we obtain the equation \eqref{8.22}. 


To estimate the coefficient of the first term on right side of \eqref{8.22}, direct computation implies
$$
\langle u_0,(r\partial_rJ)_{\lambda_0}\rangle=\left\langle \frac{u_0}{r},r(r\partial_r J)_{\lambda_0}\right\rangle
\le\left(\displaystyle\int_{\mathbb{R}^+}\left(\frac{u_0}{r}\right)^2rdr\right)^\frac{1}{2}\cdot\left(\displaystyle\int_{\mathbb{R}+}(r\partial_rJ_{\lambda_0})^2r^3dr\right)^\frac{1}{2}\le\frac{ C \v}{\lambda^2_0},$$
where we have used the assumption \eqref{2930}. 
Hence by the definition of $\la_0$ and $\dot{\lambda} (0)$, we conclude that
\beq\notag
\frac{C_0\dot\lambda(0)}{\lambda_0^2}-2\dot\lambda(0)\langle u_0,(r\partial_r J)_{\lambda_0}\rangle-2\frac{C_0}{\lambda_0}\geq C_0\v - 2C\v^{2} -2C_0/\v^{-4}\geq C\v
\eeq
for some constant $C> 0$ independent of $\v$. The estimate \eqref{sec5pf1} can be proved similarly as the estimate \eqref{sec4pf3},
which completes the proof of proposition.
\endpf

\bigskip
\noindent{\bf Proof of Theorem \ref{mthsing}}. \quad
We will prove the theorem by a contradiction argument. First assume that there is no singularity before the time $t=\v^{\frac{1}{2}}$. 
For convenience, we first rewrite the equation of \eqref{8.22} as follows 
\beq\label{8.222}
K_1\dot{\la}=K_2\la^2+2C_0\la+\la^2\int_0^t \mathcal{E}_1(s) ds
+\la^2\int_0^t \mathcal{E}_2(s) ds,
\eeq
where
$$
K_1=C_0-2\langle u,r\partial_r J\rangle, \quad 
K_2=C_0\frac{\dot\lambda(0)}{\lambda_0^2}-2\dot\lambda_0\langle u_0,(r\partial_r J)_{\lambda_0}\rangle-\frac{2C_0}{\lambda_0},
$$
\[
\begin{split}
\mathcal{E}_1=-2\langle u,\dot{J}_\lambda\rangle\dot{\lambda}-\langle u,\ddot{J}_\lambda\rangle\lambda+\langle \mathcal{N}(u),J_\lambda\rangle\lambda,\quad 
\mathcal{E}_2=-\langle u_t,J_\lambda\rangle\lambda-\langle h_t,J_\lambda\rangle\lambda.
\end{split}
\]
To obtain the contradiction, we denote  
\[\gamma(t)=-\frac{\dot{\lambda}}{\la^2}(t),\]
then the Ricciti equation \eqref{8.222} gives that 
\[
-K_1\dot{\gamma}=2C_0\gamma+\mathcal{E}_1+\mathcal{E}_2
\]

By the estimate \eqref{sec5pf1} in Lemma \ref{lemma6.1}, we know that for sufficiently small $\v$ it holds
\beq\label{thm3.1pf1}
0<C_0-C\v\leq K_1\leq C_0+C\v, 
\eeq
The terms in $\mathcal{E}_1$ are exactly the same as those in \cite{RS}. By the Proposition 5.1 and Lemma 5.6 in \cite{RS}, we have the following estimate
\beq\label{estE1}
|\mathcal{E}_1(t)|\leq C c_0^{\frac{1}{2}}
\sup_{0\leq s\leq t}\frac{\dot{\la}^4}{\la^7}(s)
=
Cc_0^{\frac{1}{2}}\frac{\dot{\la}^4}{\la^7}(t),
\eeq
since $\frac{\dot{\la}^4}{\la^7}(t)$ is monotonic increasing.
Finally similar as the proof of Lemma \ref{lemma5.5}, we have 
\[
\left|\langle u_t,J_\lambda\rangle\lambda\right|
\leq
\la\left(\int u_t^2\, rdr\right)^\frac{1}{2}\left(\int J_r^2\, rdr\right)^\frac{1}{2}
\leq C\v \leq  Cc_0^{\frac{1}{2}} 
\frac{\dot{\la}^4}{\la^7}(t),
\]
and for $t\in[0,\v^{\frac{1}{2}}]$
\[
\left|\langle h_t,J_\lambda\rangle\lambda\right|
\leq
\la\left(\int h_t^2\, rdr\right)^\frac{1}{2}\left(\int J_r^2\, rdr\right)^\frac{1}{2}
\leq C \v^{\frac{1}{4}} \leq Cc_0^{\frac{1}{2}} 
\frac{\dot{\la}^4}{\la^7}(t)\]
when $\v\leq c_0^2$ sufficiently small. Combining these two estimates, we obtain
\beq\label{estE2}
|\mathcal{E}_2(t)|\leq 
Cc_0^{\frac{1}{2}}\frac{\dot{\la}^4}{\la^7}(t).
\eeq
Therefore, choosing $\v$ small enough, and combining estimates \eqref{thm3.1pf1}-\eqref{estE2}, we obtain 
\[
\dot{\gamma} \leq -C_2\gamma +
C_2c_0^{\frac{1}{2}}\frac{\dot{\la}^4}{\la^7}
=-C_2\gamma -
C_2c_0^{\frac{1}{2}}\frac{\dot{\la}}{\la}\gamma^3,
\]
for any time $t\in[0,\v^{\frac{1}{2}}]$. Denote $\xi=e^{C_2t}\gamma$, which satisfies
\beq\notag
\dot{\xi}\leq 
-C_2c_0^{\frac{1}{2}}e^{-2C_2t}\frac{\dot{\la}}{\la}\xi^3.
\eeq
Since $-\xi>0$, dividing both sides of this inequality by $-\xi^3$ and integrating over $[0,t]$ for any $t\in[0,\v^{\frac{1}{2}}]$, it holds
\beq\notag
\frac{1}{\xi^2}\leq \frac{1}{\xi^2(0)}+Cc_0^{\frac{1}{2}}\int_0^te^{-2C_2s}\frac{\dot{\la}}{\la}(s)\,ds
\leq C\v^{-2}+Cc_0^{\frac{1}{2}}\int_0^t\frac{\dot{\la}}{\la}(s)\,ds.
\eeq
By the assumption $\la_0=\v^{-4}$, we have
\beq\notag
\int_0^t\frac{\dot{\la}}{\la}(s)\,ds=\ln(\la)-\ln(\la_0)=\ln(\la)+4\ln(\v)
\leq \ln(\la)+C\v^{-2}.
\eeq
Hence
\beq\notag
\frac{1}{\xi^2}
\leq C\v^{-2}+Cc_0^{\frac{1}{2}}\ln(\la),
\eeq
or equivalently
\[
\frac{1}{\sqrt{\v^{-2}+c_0^{\frac{1}{2}}\ln(\la)}}\leq C \frac{\dot{\lambda}}{\la^2}(t).
\]
Using the fact $\la(t)\geq\la_0= \v^{-4}$, we have
\[
\sqrt{\v^{-2}+c_0^{\frac{1}{2}}\ln(\la)}\leq C\sqrt{\la},
\]
which implies that
\[
\lambda^{\frac{3}{2}}\leq C  \dot{\lambda}(t)\
\]
when $\v$ is small enough. Since the initial data $\la_0=\v^{-4}$, solving this inequality, we can conclude that there will be a blowup  at a time $O(\v^2)$, which is a contradiction. 

The asymptotic estimates \eqref{33} can be done as in \cite{RS} (cf. page 212-213). Therefore, we complete the proof of Theorem \ref{mthsing}.



%


\bigskip
\begin{appendices}
\section{Derivation of simplified system \eqref{simeqn}}
\setcounter{equation}{0}

First note that for any function $f(x,y,z)=f(r,\theta)$
\beq\label{difff}
\begin{split}
f_x=f_r\cos \theta-f_{\theta}\frac{\sin \theta}{r},\quad
f_y=f_r\sin \theta+f_{\theta}\frac{\cos \theta}{r},\quad
f_z=0.
\end{split}
\eeq
Since
$\u=\big(0,0, v\big)^T$,
it holds
\beq\notag
\begin{split}
2D=\nabla \u+(\nabla \u)^T
=
\left[
\begin{array}{ccc}
0&0&v_r\cos k\theta\\
\\
0&0&v_r\sin k\theta\\
\\
v_r\cos k\theta&v_r\sin k\theta&0
\end{array}
\right]
\end{split}
\eeq
and
\beq\notag
\begin{split}
2\omega=\nabla \u-(\nabla \u)^T
=
\left[
\begin{array}{ccc}
0&0&-v_r\cos k\theta\\
\\
0&0&-v_r\sin k\theta\\
\\
v_r\cos k\theta&v_r\sin k\theta&0
\end{array}
\right].
\end{split}
\eeq
Thus
\beq\notag
-2\omega\d
=v_r\big(\cos\phi\cos k\theta,\cos\phi\sin k\theta,
-\sin\phi \big)^T.
\eeq

We first compute the equation of $\d=\big(\sin\phi\cos k\theta, \sin \phi\sin k\theta, \cos\phi\big)^T$. Direct calculation implies
\beq\notag
\d_t=\phi_t\big(\cos\phi\cos k\theta, \cos\phi\sin k\theta, -\sin\phi\big)^T
\eeq
\beq\notag
\d_r=\phi_r\big(\cos\phi\cos k\theta, \cos\phi\sin k\theta, -\sin\phi\big)^T
\eeq
\beq\notag
\d_{\theta}=k\big(-\sin\phi\sin k\theta, \sin\phi\cos k\theta, 0\big)^T
\eeq
and
\beq\notag
\d_{tt}=\phi_{tt}\big(\cos\phi\cos k\theta, \cos\phi\sin k\theta, -\sin\phi\big)^T+|\phi_t|^2\big(-\sin\phi\cos k\theta, -\sin\phi\sin k\theta, -\cos\phi\big)^T
\eeq
\beq\notag
\d_{rr}=\phi_{rr}\big(\cos\phi\cos k\theta, \cos\phi\sin k\theta, -\sin\phi\big)^T+|\phi_r|^2\big(-\sin\phi\cos k\theta, -\sin\phi\sin k\theta, -\cos\phi\big)^T
\eeq
\beq\notag
\d_{\theta\theta}=k^2\big(-\sin\phi\cos k\theta, -\sin\phi\sin k\theta, 0\big)^T
\eeq
Hence
\beq\notag
\u\cdot\nabla \d
=0,
\eeq
\beq\notag
\begin{split}
\Delta \d=&\d_{rr}+\frac{\d_r}{r}+\frac{\d_{\theta\theta}}{r^2}\\
=&\left(\phi_{rr}+\frac{\phi_r}{r}\right)
\big(\cos\phi\cos k\theta, \cos\phi\sin k\theta, -\sin\phi\big)^T\\
&-|\phi_r|^2\big(\sin\phi\cos k\theta, \sin\phi\sin k\theta, \cos\phi\big)^T\\
&-\frac{k^2}{r^2}\big(\sin\phi\cos k\theta, \sin\phi\sin k\theta, 0\big)^T,
\end{split}
\eeq
and
\beq\notag
\begin{split}
&\big(|\nabla \d|^2-|\dot{\d}|^2\big)\d\\
=&\left(|\d_r|^2+\frac{|\d_{\theta}|^2}{r^2}-|\d_t|^2\right)\d\\
=&\left(|\phi_r|^2+\frac{\sin^2\phi}{r^2}-|\d_t|^2\right)\d\\
=&(|\phi_r|^2-|d_t|^2)\big(\sin\phi\cos k\theta, \sin\phi\sin k\theta, \cos\phi\big)^T\\
&+\frac{k^2\sin^2\phi}{r^2}\big(\sin\phi\cos k\theta, \sin\phi\sin k\theta, \cos\phi\big)^T.
\end{split}
\eeq
It is easy to see that
\beq\notag
\begin{split}
&\ddot{\d}-\Delta \d-\big(|\nabla \d|^2-|\dot{\d}|^2\big)\d\\
=&\left(\phi_{tt}-\phi_{rr}-\frac{\phi_r}{r}\right)
\big(\cos\phi\cos k\theta, \cos\phi\sin k\theta, -\sin\phi\big)^T\\
&+\frac{k^2}{r^2}\big(\sin\phi\cos k\theta, \sin\phi\sin k\theta, 0\big)^T-\frac{k^2\sin^2\phi}{r^2}\big(\sin\phi\cos k\theta, \sin\phi\sin k\theta, \cos\phi\big)^T\\
=&\left(\phi_{tt}-\phi_{rr}-\frac{\phi_r}{r}+\frac{k^2\sin\phi\cos\phi}{r^2}\right)
\big(\cos\phi\cos k\theta, \cos\phi\sin k\theta, -\sin\phi\big)^T.
\end{split}
\eeq
On the other hand
\beq\notag
\begin{split}
2\dot \d-2\omega \d
=(2\phi_t+v_r)
\big(\cos\phi\cos\theta, \cos\phi\sin\theta, -\sin\phi\big)^T.
\end{split}
\eeq
Therefore the equation of $d $ in system \eqref{fels} becomes
\beq
\phi_{tt}+2\phi_t=\phi_{rr}+\frac{\phi_r}{r}-\frac{k^2\sin\phi\cos\phi}{r^2}-v_r.
\eeq

Now, we proceed to work on the equation of $\u$. Since
$\u=\big(0, 0, v\big)^T$, direct calculation implies
\beq\notag
\u_t=\big(0,0, v_t\big)^T,\quad\u_r=\big(0,0, v_r\big)^T,\quad \u_{rr}=\big(0,0, v_{rr}\big)^T,\quad\u_{\theta}=\u_{\theta\theta}=\big(0,0, 0\big)^T.
\eeq
Hence
\beq\notag
\nabla\cdot \u=\u\cdot\nabla \u=0,\quad\Delta \u=\u_{rr}+\frac{\u_r}{r}+\frac{\u_{\theta\theta}}{r^2}
=
\left(0,0,
v_{rr}+\frac{v_r}{r}\right)^T.
\eeq
For the other terms, first notice that
\beq\notag
\nabla P
=
\big(P_r\cos k\theta, P_r\sin k\theta,0\big)^T.
\eeq
By the definition of $N$, we have
\beq\notag
\begin{split}
N=\d_t+\u\cdot\nabla \d-\omega \d
=\left(\phi_t+\frac12v_r\right)
\big(\cos\phi\cos k\theta, \cos\phi\sin k\theta, -\sin\phi\big)^T.
\end{split}
\eeq
Thus
\beq\notag
\begin{split}
\d\otimes N-N\otimes \d
=\left(\phi_t+\frac12v_r\right)
\left[
\begin{array}{ccc}
0&0&-\cos k\theta\\
\\
0&0&-\sin k\theta\\
\\
\cos k\theta&\sin k\theta&0
\end{array}
\right],
\end{split}
\eeq
and furthermore
\beq\notag
\begin{split}
\nabla\cdot\big(\d\otimes N-N\otimes \d\big)
=
\left(
0,0,
\displaystyle \frac1r\left[r\left(\phi_t+\frac12v_r\right)\right]_r
\right)^T.
\end{split}
\eeq
By the definition of $D$ and the fact $\nabla\cdot \u=0$, we obtain
\beq\notag
\begin{split}
&2\nabla \cdot D=\partial_{x_j}\big(\u^i_{x_j}+\u^j_{x_i}\big)=\Delta u
=\left(0,0,
v_{rr}+\frac{v_r}{r}\right)^T.
\end{split}
\eeq
For the last term in the equation of $\u$, it can be computed as follows
\beq\notag
\begin{split}
&\nabla \cdot (\nabla \d\odot \nabla \d)=\partial_{x_j}\big(\d_{x_i}\cdot \d_{x_j}\big)=\nabla \d\cdot \Delta \d+\frac12\nabla\big(|\nabla \d|^2\big).
\end{split}
\eeq
It is not hard to see
\beq\notag
\nabla\big(|\nabla \d|^2\big)=\nabla\left(|\phi_r|^2+\frac{\sin^2\phi}{r^2}\right)=\left(\big(|\nabla \d|^2\big)_r\cos k\theta,
\big(|\nabla \d|^2\big)_r\sin k\theta, 0 \right)^T
\eeq
and
\beq\notag
\begin{split}
\nabla \d\cdot \Delta \d
=\left(
\begin{array}{ccc}
\displaystyle\left(\phi_{rr}+\frac{\phi_r}{r}-\frac{k^2\sin\phi\cos\phi}{r^2}\right)\phi_r\cos k\theta,
 \displaystyle\left(\phi_{rr}+\frac{\phi_r}{r}-\frac{k^2\sin\phi\cos\phi}{r^2}\right)\phi_r\sin k\theta,
 0
\end{array}
\right)^T.
\end{split}
\eeq
So the equation of $v$ is given by
\beq
\begin{split}
v_t=\frac{1}{r}\Big(rv_r+r\phi_t\Big)_r.
\end{split}
\eeq



\section{A stronger decay estimate}
\setcounter{equation}{0}

This section is devoted to prove the stronger decay estimate for the solutions to \eqref{simeqn}.

\begin{lemma}\label{lemma6.4}
Assume that $(\phi, v)$ exists and is smooth on time interval $(0,T]$. Let $\lambda(t)$ be the solution obtained in Proposition \ref{prop3.1} with additional assumption of $\lambda\geq 11$. Denote $u=\phi-I_\lambda$ which satisfies the smallness condition \eqref{cond1}. Then the following stronger estimate of $u$ is valid for any time $t\in(0,T]$
\beq\label{148}
\int_{r\geq 2t} r^2\left((\phi_t)^2+(u_r)^2+\frac{k^2}{r^2}u^2+v^2\right)\,rdr\leq C \v^2.
\eeq 
\end{lemma}

\pf  It is not hard to see that we only need to show \eqref{148} in the region $\{r\geq 2t+11\}$. Denote $\alpha(y)$ as a smooth increasing function supported on $y\geq 10$ and satisfying $\alpha(y)=y^2$ when $y\geq11$ and
\beq
\alpha'(y)\leq 3y,\quad y^{-1}\alpha(y)\leq \alpha'(y),\quad |\alpha''|\leq C.
\eeq
By the second equation of \eqref{simeqn} and the notation in \eqref{3.39}, we have
\beq\label{149}
\int_0^t\int_{\R^+}\left(\phi_{tt}+2\phi_t+H_{\lambda}u-\mathcal N(u)+v_r\right)\phi_t\,\alpha(r-2s)\,rdrds=0.
\eeq
By the definition of the Hamiltonian, we denote 
\beq\notag
H_\lambda=-\partial_{rr}-\frac{\partial_r}{r}+\frac{k^2}{r^2}\cos(2I_\lambda):=-\partial_{rr}-\frac{\partial_r}{r}+Q_{\lambda}.
\eeq
It is not hard to see that $Q_\lambda\geq C\frac{k^2}{r^2}$ on the support of $\alpha(r-2s)$.

For the first term in \eqref{149}, integrating by parts implies 
\beq\label{149-1}
\int_0^t\int_{\R^+}\phi_{tt}\phi_t\,\alpha(r-2s)\,rdrds
=\frac12\int_{\R^+}|\phi_t|^2\alpha(r-2s)\,rdr\big|_{s=0}^{s=t}
+\int_0^t\int_{\R^+}|\phi_t|^2\alpha'(r-2s)\,rdrds.
\eeq
Using $\phi=u+I_\lambda$, we may rewrite the third term in \eqref{149} as follows
\beq\notag
\int_0^t\int_{\R^+}H_{\lambda}u\phi_t\,\alpha(r-2s)\,rdrds
=\int_0^t\int_{\R^+}(H_{\lambda}u) u_t\,\alpha(r-2s)\,rdrds
+\int_0^t\int_{\R^+}(H_{\lambda}u) (I_\lambda)_t\,\alpha(r-2s)\,rdrds.
\eeq
Integrating by parts implies
\beq\notag
\begin{split}
&\int_0^t\int_{\R^+}(H_{\lambda}u) u_t\,\alpha(r-2s)\,rdrds\\
=&\int_0^t\int_{\R^+}(u_r u_{rt}\,\alpha(r-2s)+u_ru_t\alpha'(r-2s))\,rdrds+\int_0^t\int_{\R^+}(Q_\lambda u)u_t\alpha(r-2s)\,rdrds\\
=&\frac12\int_{\R^+}\left(|u_r|^2+Q_\lambda u^2\right)\alpha(r-2s)\,rdr\big|_{s=0}^{s=t}
+\int_0^t\int_{\R^+}\left(|u_r|^2+u_ru_t+Q_\lambda u^2\right)\alpha'(r-2s)\,rdrds\\
&-\frac12\int_0^t\int_{\R^+}\dot{Q}_\lambda u^2\alpha(r-2s)\,rdrds.
\end{split}
\eeq
By the decomposition $H_\lambda=A^*_\lambda A_\lambda$
and the identity \eqref{AId1}, we have
\beq\notag
\int_0^t\int_{\R^+}(H_{\lambda}u) (I_\lambda)_t\,\alpha(r-2s)\,rdrds
=\int_0^t\int_{\R^+}\left(u_r-\frac{ku}{r}\cos(I_\lambda)\right)(I_\lambda)_t\alpha'(r-2s)\,rdrds.
\eeq
Hence the third term in \eqref{149} can be estimated as follows
\beq\label{149-2}
\begin{split}
&\int_0^t\int_{\R^+}H_{\lambda}u\phi_t\,\alpha(r-2s)\,rdrds\\
=&\frac12\int_{\R^+}\left(|u_r|^2+Q_\lambda u^2\right)\alpha(r-2s)\,rdr\big|_{s=0}^{s=t}
+\int_0^t\int_{\R^+}\left(|u_r|^2+u_r\phi_t+Q_\lambda u^2\right)\alpha'(r-2s)\,rdrds\\
&-\frac12\int_0^t\int_{\R^+}\dot{Q}_\lambda u^2\alpha(r-2s)\,rdrds-\int_0^t\int_{\R^+}\frac{ku}{r}\cos(I_\lambda)(I_\lambda)_t\alpha'(r-2s)\,rdrds.
\end{split}
\eeq
Combining all these estimates with \eqref{149}, we obtain
\beq\label{150-1}
\begin{split}
&\frac12\int_{\R^+}\left(|\phi_t|^2+|u_r|^2+Q_\lambda u^2\right)\alpha(r-2s)\,rdr\big|_{s=0}^{s=t}\\
=&-\int_0^t\int_{\R^+}\left(|\phi_t|^2+|u_r|^2+u_r\phi_t+Q_\lambda u^2\right)\alpha'(r-2s)\,rdrds\\
&+\frac12\int_0^t\int_{\R^+}\dot{Q}_\lambda u^2\alpha(r-2s)\,rdrds+\int_0^t\int_{\R^+}\frac{ku}{r}\cos(I_\lambda)(I_\lambda)_t\alpha'(r-2s)\,rdrds\\
&-\int_0^t\int_{\R^+}\left(2|\phi_t|^2+v_r\phi_t\right)\alpha(r-2s)\,rdrds+\int_0^t\int_{\R^+}\mathcal{N}(u)\phi_t\alpha(r-2s)\,rdrds
\end{split}
\eeq
By the first equation of \eqref{simeqn}, we have
\beq\notag
\int_0^t\int_{\R^+}\left(v_t-\frac{1}{r}\Big(rv_r+r\phi_t\Big)_r\right)v\,\alpha(r-2s)\,rdrds=0.
\eeq
Integrating by parts implies
\beq\label{150-2}
\begin{split}
&\frac12\int_{\R^+}|v|^2\alpha(r-2s)\,rdr\big|_{s=0}^{s=t}\\
=&-\int_0^t\int_{\R^+}\left(|v_r|^2+\phi_tv_r\right)\alpha(r-2s)\,rdrds-\int_0^t\int_{\R^+}\left(|v|^2+\phi_tv\right)\alpha'(r-2s)\,rdrds\\
&-\int_0^t\int_{\R^+}v_rv\alpha'(r-2s)\,rdrds.
\end{split}
\eeq
Combining \eqref{150-1} and \eqref{150-2}, we obtain
 \beq\label{150}
\begin{split}
&\frac12\int_{\R^+}\left(|\phi_t|^2+|u_r|^2+Q_\lambda u^2+|v|^2\right)\alpha(r-2s)\,rdr\big|_{s=0}^{s=t}\\
\leq&-\int_0^t\int_{\R^+}\left(|\phi_t|^2+|u_r|^2+|v_r|^2+\phi_tv_r+u_r\phi_t+\frac{Ck^2}{r^2} u^2\right)\alpha'(r-2s)\,rdrds\\
&-\int_0^t\int_{\R^+}\left(2|\phi_t|^2+2v_r\phi_t+|v_r|^2\right)\alpha(r-2s)\,rdrds\\
&+\frac12\int_0^t\int_{\R^+}\dot{Q}_\lambda u^2\alpha(r-2s)\,rdrds+\int_0^t\int_{\R^+}\mathcal{N}(u)\phi_t\alpha(r-2s)\,rdrds\\
&+\int_0^t\int_{\R^+}\frac{ku}{r}\cos(I_\lambda)(I_\lambda)_t\alpha'(r-2s)\,rdrds
-\int_0^t\int_{\R^+}v_rv\alpha'(r-2s)\,rdrds\\
:=&\sum\limits_{j=1}^6T_j,
\end{split}
\eeq
where we have use the fact that $Q_\lambda\geq C\frac{k^2}{r^2}$ on the support of $\alpha(r-2s)$.

It is not hard to see that $T_1\geq 0$ and $T_2\geq 0$. By the definition of $Q_\lambda$ and \eqref{cond3} in Proposition \ref{prop3.1}, we have
\beq\notag
|\dot{Q}_{\lambda}|\leq C \frac{\v}{r^3}.
\eeq
Combining the property of $\alpha$ implies
\beq\notag
\left|T_3\right|\leq C \int_0^t\int_{\{r-2s\geq 10\}} \frac{\v}{r^3}u^2(r-2s)\alpha'(r-2s)\,rdrds\leq C \v T_1.
\eeq
By the definition of $\mathcal{N}(u)$, inequality \eqref{68} and \eqref{cond1} in Proposition \ref{prop3.1}, we can estimate $T_4$ as follows
\beq\notag
\left|T_4\right|\leq C \int_0^t\int_{\{r-2s\geq 10\}} \v\frac{|u|}{r^2}|\phi_t|(r-2s)\alpha'(r-2s)\,rdrds
\leq C \v T_1.
\eeq
Therefore, for sufficiently small $\v>0$, we have
\beq
T_1+T_3+T_4\geq (1-C\v)T_1\geq 0.
\eeq
For term $T_5$, by the definition of $I_\lambda$, the assumptions $k\geq 4$ and $\lambda\geq 1$, and the estimate \eqref{cond3} in Proposition \ref{prop3.1}, we obtain
\beq\notag
|\dot{I}_\lambda \alpha'|\leq C|r\dot{I}|\leq C\frac{\v}{(1+r)^3}.
\eeq
Hence, we can estimate $T_5$ as follows
\beq
\begin{split}
\left|T_5\right|
\leq C \left(\int_0^t\int_{\R^+} \frac{|u|^2}{r^2}\,rdrds\right)^{\frac12}\left(\int_0^t\int_{\R^+} |\dot{I}_\lambda|^2\alpha'(r-2s)\,rdrds\right)^{\frac12}
\leq C \v^2, 
\end{split}
\eeq
where we has used the estimate \eqref{cond1} in Proposition \ref{prop3.1} in last inequality. For term $T_6$, integrating by parts implies
\beq\notag
\begin{split}
T_6=\frac12\int_0^t\int_{\R^+}|v|^2\alpha''(r-2s)\,rdrds
+\frac12\int_0^t\int_{\R^+}|v|^2\alpha'(r-2s)\,drds.
\end{split}
\eeq
Hence by the property of $\alpha$ and the estimate \eqref{cond1} in Proposition \ref{prop3.1}
\beq
\begin{split}
|T_6|\leq C\int_0^t\int_{\R^+}|v|^2\,rdrds\leq C \v^2.
\end{split}
\eeq
Putting all these estimates in \eqref{150}, we can conclude the estimate \eqref{148}, which completes the proof.

\endpf

We are ready to show a stronger pointwise decay estimate of $u$ outside of a sufficiently large cone. 

\begin{proposition}\label{decayu}
Suppose that $u$ satisfies the assumption of Lemma \ref{lemma6.4}. The following pointwise estimates are valid 
\beq\label{147}
\sup\limits_{r,t}|u(r,t)|\leq C\v,\quad \sup\limits_{r\geq 3t}r^2|u|^2(r,t)\leq C \v^2
\eeq
for all $t\in [0,T]$. 
\end{proposition}

\pf First notice that when $t=0$, the estimate \eqref{147}  can be deduced by the assumptions on initial data $u_0$. Hence, we only need to investigate $t\in(0,T]$.
By the Poincar\'e type estimate \eqref{68} and \eqref{cond1}, it is easy to see that 
\beq\label{infbu}
(u(r))^2
\le 2\left(\displaystyle\int_{\mathbb{R}^+}\left(\frac{u(r)}{r}\right)^2rdr\right)^\frac{1}{2}\left(\displaystyle\int_{\mathbb{R}^+}(u_r)^2rdr\right)^\frac{1}{2}\leq 2\v^2,
\eeq
which implies the first estimate in \eqref{147}.

Denote $\chi(r,t)$ is a smooth cutoff function on the region $r\geq 2t$ which satisfies $\chi\equiv 1$ on $r\geq 3t$, $\chi\equiv 0$ outside of $r\geq 2t$, $|\chi|\leq C$ and $|\chi_r|\lesssim r^{-1}$. Applying the Poincar\'e type estimate \eqref{68} to $r\chi u(r)$, we have 
\beq\label{infbur}
\begin{split}
&\sup\limits_{r\geq 3t}r^2|u|^2\leq \sup\limits_{r\geq 3t} |r\chi u(r)|^2\\
\leq & C\left(\int_{\mathbb{R}^+}|\chi u|^2rdr\right)^\frac{1}{2}\left(\displaystyle\int_{\mathbb{R}^+}(|r\chi u_r|^2+|\chi u|^2+|r\chi_r u|^2)rdr\right)^\frac{1}{2}\\
\leq & C \left(\int_{\mathbb{R}^+}|u|^2rdr\right)^\frac{1}{2}\left(\displaystyle\int_{r\geq 2t}(r^2|u_r|^2+|u|^2)rdr\right)^\frac{1}{2}\\
\leq & C\v^2,
\end{split}
\eeq
where we have used \eqref{148} in last inequality. This completes the proof of \eqref{147}.


\end{appendices}

\section*{Acknowledgments}
The first and third authors are partially supported by NSF with grant DMS-2008504.

\end{document}